\documentclass[12pt]{amsart}
\usepackage{amsmath,amsfonts,amsthm,mathrsfs,mathtools,amscd}
\usepackage{amssymb,latexsym}
\usepackage{cases}
\usepackage[numbers]{natbib}
\usepackage[shortlabels]{enumitem}
\setlist[enumerate]{nosep}
\usepackage{hyperref}
\usepackage{cleveref}
\newtheorem{thm}{Theorem}[section]

\newtheorem{lem}[thm]{Lemma}
\newtheorem{prop}[thm]{Proposition}
\newtheorem{cor}[thm]{Corollary}
\theoremstyle{remark}
\newtheorem{re}[thm]{Remark}
\theoremstyle{definition}

\newtheorem{defn}[thm]{Definition}

\newcommand{\cls}{\mbox{cls}}
\newcommand{\spn}{\mbox{spn}}
\newcommand{\gen}{\mbox{gen}}
\newcommand{\ord}{\mathrm{ord}}

\newcommand{\rank}{\mbox{rank}}

\newcommand{\tr}{\mbox{Tr}}
\newcommand{\nr}{\mbox{N}}

\newcommand{\gal}{\mbox{Gal}}
\newcommand{\id}{\mbox{Id}}
\newcommand{\cl}{\mbox{Cl}}

\numberwithin{equation}{section}

\ExplSyntaxOn
\seq_new:N \l__ne_bib_entries_seq
\cs_new:Npn \__ne_bib_key:w [#1]#2\q_stop {#1}
\cs_new_protected:Npn \__ne_bib_print:w [#1]#2\q_stop
  { \bibitem #2 }
\NewDocumentEnvironment{sortedbibliography}{m +b}
  {
    \seq_set_split:Nnn \l__ne_bib_entries_seq {\bibitem} {#2}
    \seq_pop_left:NN \l__ne_bib_entries_seq \l_tmpa_tl
    \seq_sort:Nn \l__ne_bib_entries_seq
      {
        \str_compare:eNeTF
          { \__ne_bib_key:w ##1\q_stop } >
          { \__ne_bib_key:w ##2\q_stop }
          { \sort_return_swapped: }
          { \sort_return_same: }
      }
    
  }
  {}
\ExplSyntaxOff

\begin{document}
\title[]{Class numbers and finiteness of norm Euclidean quadratic lattices over algebraic number fields}
\author{Zilong He}
\address{Department of Mathematics,
	Dongguan University of Technology, Dongguan 523808, China}
\email{zilonghe@connect.hku.hk}

\subjclass[2020]{Primary 11E12, 11E41; Secondary 11R29}
\date{}
\keywords{norm Euclidean lattices, class number one, maximal lattices}

\begin{abstract}
	Let $F$ be an algebraic number field and $\mathcal{O}_{F}$ the ring of integers of $F$. Let $L$ be an integral
	$\mathcal{O}_{F}$-lattice of rank $n\ge1$ on an anisotropic quadratic space $V$.
	We prove that if $L$ is norm Euclidean, then $L$ is $\mathcal{O}_{F}$-maximal of class number one, and the class number of $F$ equals the class number of $L\perp\mathbf{H}$, where
	$\mathbf{H}$ is the hyperbolic lattice with quadratic form $xy$.  Under the same assumption, $F$ has class number one if $n\ge3$ and at most two if $n=1,2$. We also prove finiteness up to isometry for totally definite norm Euclidean lattices over totally real number
	fields in ranks at
	least three as the field varies, in rank two for fixed field
	degree, and in rank one over a fixed field.
	This confirms Clark's finiteness conjecture in ranks at least three. 
\end{abstract}

\maketitle

\section{Introduction}
\label{sec:preliminaries}

Let $F$ be a totally real number field of degree $d$ with ring of
integers $\mathcal{O}_{F}$. A totally positive definite integral
quadratic form $Q$ in $n$ variables over $\mathcal{O}_{F}$ is called
\emph{norm Euclidean}, or simply \emph{Euclidean}, if for every
$x\in F^{n}$, there exists $y\in\mathcal{O}_{F}^{n}$ such that
\begin{align*}
	\nr_{F/\mathbb{Q}}(Q(x-y))<1.
\end{align*}
This notion was first introduced and studied by Clark
\cite{clark_ADC-I-2012}.
Clark conjectured that there are only finitely many isometry classes
of totally definite integral norm Euclidean quadratic forms over
$\mathcal{O}_{F}$ as $F$ varies over all totally real number fields
\cite[Conjecture 29]{clark_ADC-I-2012}.
For $F=\mathbb{Q}$, Clark and Jagy proved that every positive
definite integral norm Euclidean quadratic form has class number
one \cite[Theorem 6.1]{clark_ADC-II-2014}. Actually,
there are exactly $70$ isometry classes of such lattices,
as classified by Nebe and recovered by Clark and Jagy;
see \cite[Theorem 1.6, Remark 1.7, and \S6.4]{clark_ADC-II-2014}.

In \cite[Definition 1.1]{he_trace_2026}, we introduced
the trace Euclidean property by replacing the norm inequality
above with the trace inequality
\begin{align*}
	\tr_{F/\mathbb{Q}}(Q(x-y))<d.
\end{align*}
This property is stronger than the norm Euclidean property
\cite[(1.2)]{he_trace_2026}.
We obtained bounds for the norms of volume ideals by relating
covolumes to trace covering radii, and used these bounds to
prove finiteness for positive definite trace Euclidean lattices
in both the classic integral and integral cases
\cite[Lemma 3.4 and Theorems 1.2 and 1.3]{he_trace_2026}. However, a norm Euclidean lattice need not be trace Euclidean. For example, over $F=\mathbb{Q}(\sqrt{3})$, the rank one lattice
$L=\mathcal{O}_{F}v$ with $Q(v)=1$
is norm Euclidean but not trace Euclidean
(cf. \cite[Theorem 1.7 and the paragraph below]{he_trace_2026}).
The finiteness theorems proved there require the trace Euclidean
hypothesis and therefore do not apply directly to norm Euclidean
lattices in general.

In this paper, we study norm Euclidean lattices on anisotropic
quadratic spaces over algebraic number fields
(see Definition \ref{defn:norm-euclidean}). Following the work of Clark and Jagy \cite[\S6.3]{clark_ADC-II-2014},
we develop the Lorentzian lattice method in this setting.
We prove maximality and class number one for such lattices
and obtain restrictions on the class number of the base field
(Theorem \ref{thm:class-number}). In particular, this proves Clark's class number conjecture
\cite[Conjecture 28]{clark_ADC-I-2012} for rings of integers
of number fields equipped with their canonical norms, and recovers
\cite[Theorem 6.1]{clark_ADC-II-2014} when $F=\mathbb{Q}$ and $L$ is positive definite.
A key ingredient in our proof is the cancellation theorem
(Theorem \ref{thm:cancellation}).
We also prove a version with an ideal parameter
(Corollary \ref{cor:ideal-cancellation}).
These results are analogous to the cancellation and
classification theorems of Schulze-Pillot
for unimodular lattices over Hasse domains (cf. 
\cite[Theorems 1 and 2]{schulze-pillot_remark_1982}). 

We then prove finiteness results for totally definite norm
Euclidean lattices over totally real fields
(Theorem \ref{thm:finiteness}).
This confirms Clark's finiteness conjecture in ranks at least three.
When $F$ is totally real and $L$ is a positive definite
norm Euclidean $\mathcal{O}_{F}$-lattice of rank $n\ge3$,
we also obtain explicit bounds on the rank
and the field degree (Theorem \ref{thm:rank-degree-bounds}).

Before stating the main results, we review the basic definitions
and notation. Let $F$ be an algebraic number field, and let $\mathcal{O}_{F}$ be the ring
of integers of $F$. Let $d=[F:\mathbb{Q}]$, and denote by
$\nr_{F/\mathbb{Q}}$ the field norm from $F$ to $\mathbb{Q}$. As usual, we write
$\nr(\mathfrak{a})=|\mathcal{O}_{F}/\mathfrak{a}|$ for the absolute
norm of a nonzero integral ideal $\mathfrak{a}$, and for a nonzero fractional ideal $\mathfrak{b}$ of $F$, we choose
$\gamma\in F^{\times}$ such that $\gamma\mathfrak{b}$ is integral, and put
\begin{align*}
	\nr(\mathfrak{b})
	\coloneq
	\dfrac{\nr(\gamma\mathfrak{b})}{|\nr_{F/\mathbb{Q}}(\gamma)|}.
\end{align*}
This definition is independent of the choice of $\gamma$;
see \cite[\S9.4, pp. 231--232]{alaca_introductory_2004}.
Let $\cl(F)$ be the ordinary ideal class group and put
$h(F)=|\cl(F)|$. Let $H_{F}$ be the Hilbert class field corresponding
to $\cl(F)$, in which all infinite primes of $F$ split completely.

Let $V$ be a non-degenerate $n$-dimensional quadratic space over $F$
associated with symmetric bilinear form $B$ and quadratic form $Q$
such that
\begin{align*}
	B(x,y)=\dfrac{1}{2}(Q(x+y)-Q(x)-Q(y)),
	\qquad x,y\in V.
\end{align*}
We write $\det V\in F^{\times}/F^{\times2}$ for the determinant of $V$.
Let $L$ be an $\mathcal{O}_{F}$-lattice on $V$, and write $FL$ for the subspace of $V$ spanned by $L$ over $F$.

As in \cite[\S82E]{omeara_quadratic_1963}, we denote the scale ideal, norm ideal and volume ideal of $L$ by
$\mathfrak{s}(L)$, $\mathfrak{n}(L)$ and $\mathfrak{v}(L)$, respectively. 
We call $L$ \emph{integral} if
$\mathfrak{n}(L)\subseteq\mathcal{O}_{F}$.
An integral lattice $L$ is called \emph{$\mathcal{O}_{F}$-maximal} 
if there is no integral $\mathcal{O}_{F}$-lattice $L^{\prime}$ on $V$
such that $L\subsetneq L^{\prime}$.
If $F$ is totally real, let $\sigma_{1},\ldots,\sigma_{d}$ be the
$d$ distinct real embeddings from $F$ into $\mathbb{R}$.
We call $L$ \emph{totally definite} if $Q$ is positive or negative definite at every $\sigma_{i}$.
We call $L$ \emph{totally positive definite}, or simply \emph{positive definite},
if $Q$ is positive definite at every $\sigma_{i}$. 

We denote by $O(V)$ the orthogonal group of $V$, and put
\begin{align*}
	O(L)=\{\sigma\in O(V)\mid \sigma L=L\}
\end{align*}
if $V=FL$. Write $O^{+}(V)$ for the group of isometries of determinant $1$ on a
quadratic space $V$, and put $O^{+}(L)=O(L)\cap O^{+}(V)$. We use the notation $\cls(L)$, $\cls^{+}(L)$, $\gen(L)$,
$\spn(L)$ and $\spn^{+}(L)$ for the class, proper class,
genus, spinor genus and proper spinor genus of $L$, respectively;
see \cite[\S82C and \S102A]{omeara_quadratic_1963}.
We write $h(L)$, $h^{+}(L)$, $g(L)$ and $g^{+}(L)$ for the
numbers of classes, proper classes, spinor genera and proper
spinor genera in $\gen(L)$, respectively.

For a nonzero fractional ideal $\mathfrak{a}$ of $F$, let
\begin{align*}
	\mathbf{H}(\mathfrak{a})
	\coloneq\mathfrak{a}e+\mathfrak{a}^{-1}f 
\end{align*}
be a binary lattice on the usual
hyperbolic plane $\mathbb{H}$ with
\begin{align*}
	Q(e)=Q(f)=0
	\quad\text{and}\quad
	B(e,f)=\dfrac{1}{2}.
\end{align*} 
Thus $Q(\alpha e+\beta f)=\alpha\beta\in\mathcal{O}_{F}$
for $\alpha\in\mathfrak{a}$ and $\beta\in\mathfrak{a}^{-1}$,
so $\mathbf{H}(\mathfrak{a})$ is integral.
Put $\mathbf{H}=\mathbf{H}(\mathcal{O}_{F})
=\mathcal{O}_{F}e+\mathcal{O}_{F}f$, and
\begin{align}\label{lorentzian}
	M=L\perp\mathbf{H}.
\end{align}
In other words, $M$ splits off the hyperbolic lattice $\mathbf{H}$.
If $F$ is totally real and $L$ is positive definite of rank $n$,
then $M$ has signature $(n+1,1)$ at every real prime.
When $F=\mathbb{Q}$ and $L$ is positive definite,
$M$, equipped with the bilinear form $2B$, is the Lorentzian
lattice considered in \cite[\S6.3]{clark_ADC-II-2014};
see also \cite[\S1]{borcherds_lorentzian_1987}.

\begin{defn}\label{defn:norm-euclidean}
	Let $L$ be an integral $\mathcal{O}_{F}$-lattice on an anisotropic
	quadratic space $V$. We say that $L$ is norm Euclidean 
	if for every $x\in V$, there exists $y\in L$ such that
	\begin{align*} 
		|\nr_{F/\mathbb{Q}}(Q(x-y))|<1.
	\end{align*}
\end{defn}
For a free lattice, this agrees with Clark's definition for
anisotropic forms equipped with the canonical norm
\cite[\S 1.3 and \S 2.1]{clark_ADC-I-2012}.
\begin{thm}\label{thm:class-number}
 	Let $F$ be an algebraic number field and let $L$ be an
 	$\mathcal{O}_{F}$-lattice of rank $n\ge1$ on an anisotropic quadratic space $V$. Put $E=F(\sqrt{-\det V})$. Suppose that $L$ is norm Euclidean.
 	Then $L$ is $\mathcal{O}_{F}$-maximal, $h(L)=1$, and $h(F)=h(L\perp\mathbf{H})$. More precisely,
 	\begin{enumerate}[label=\rm(\roman*), itemindent=0em, leftmargin=*]
 		\item If $n\ge3$, then $h(F)=1$.
 		\item If $n=2$, then $h(F)\leq2$, with equality if and only if $H_{F}=E$.
 		\item If $n=1$, then $h(F)\leq2$, with equality if and only if $L$ is not free.
 	\end{enumerate}
\end{thm} 
\begin{re}\label{re:class-number}
	(i) Suppose that $n=2$. If $E\not\subseteq H_{F}$, then by
	Theorem \ref{thm:class-number}(ii), $h(F)=1$.
	If $V$ is definite at a real prime $\mathfrak{p}$, then $\mathfrak{p}$ becomes complex in $E$
	but splits completely in $H_{F}$, so $E\not\subseteq H_{F}$ and $h(F)=1$.
	The same conclusion holds if $E/F$ is ramified at a finite prime.
	
	(ii) For the standard rank one lattice $J=\mathcal{O}_{F}$
	with $Q(x)=x^{2}$, norm Euclideanity of $J$ is equivalent
	to norm Euclideanity of $F$ (cf. \cite[(1.3)]{he_trace_2026}).
	If $J$ is norm Euclidean, then $\mathcal{O}_{F}$ is a Euclidean
	domain and hence a PID, so $h(F)=1$. Theorem \ref{thm:class-number}(i) may therefore be viewed as
	a higher-dimensional analogue of this classical implication. 
\end{re}

\begin{cor}\label{cor:stabilized-genus}
	Suppose that $L$ is norm Euclidean of rank $n\ge1$ and $h(F)=2$.
	For any nonprincipal fractional ideal $\mathfrak{a}$ of $F$, we have
	\begin{align*}
		\gen(L\perp\mathbf{H})
		&=\cls(L\perp\mathbf{H})
		\cup\cls(L\perp\mathbf{H}(\mathfrak{a})).
	\end{align*}
\end{cor} 

Combining Theorem \ref{thm:class-number} with Pfeuffer's finiteness
theorems \cite[Theorems 1--3]{pfeuffer_class-numbers_1979} and
Lemma \ref{lem:bounded-volume-finiteness}, we obtain the following finiteness theorem.
\begin{thm}\label{thm:finiteness}
	The following finiteness assertions hold up to isometry.
	\begin{enumerate}[label=\rm(\roman*), itemindent=0em, leftmargin=*]
		\item\label{item:high-rank}
		As $F$ varies over all totally real number fields, there are
		only finitely many totally definite integral norm Euclidean
		$\mathcal{O}_{F}$-lattices of arbitrary rank $ \ge3$.
		
		\item\label{item:binary-bounded}
		As $F$ varies over all totally real number fields with a fixed degree $d$, there are only finitely
		many totally definite integral norm Euclidean $\mathcal{O}_{F}$-lattices
		of rank $2$.
		
		\item\label{item:fixed-field}
		For every fixed totally real number field $F$, there are only finitely
		many totally definite integral norm Euclidean $\mathcal{O}_{F}$-lattices
		of rank $1$.
	\end{enumerate}
\end{thm}
\begin{re}\label{re:unary-quadratic-finiteness}
	The finiteness assertion in Theorem \ref{thm:finiteness}(ii)
	also holds in rank one when $d=1$ or $2$.
	Indeed, in rank one, the proof of Lemma \ref{lem:field-class-group}(iii)
	shows that $F$ has a norm Euclidean ideal class.
	By \cite[(2.5), p. 129]{lenstra_euclidean_ideals_1979}, there are exactly
	$19$ real quadratic fields with such a class: $16$ have norm Euclidean
	rings of integers, and the remaining three have discriminants
	$40$, $60$ and $85$.
	Thus only finitely many fields occur when $d\le2$, and
	Theorem \ref{thm:finiteness}(iii) shows finiteness over each of them.
\end{re}
\begin{thm}\label{thm:rank-degree-bounds}
 	Let $F$ be a totally real number field of degree $d$ and let $L$ be a positive definite
 	$\mathcal{O}_{F}$-lattice of rank $n\ge3$. If $L$ is norm Euclidean, then $d\le 5$ and
 		\begin{align*}
 			n\le 
 			\begin{cases}
 				10,&\text{if $F=\mathbb{Q}$},\\
 				6,&\text{if $F\neq\mathbb{Q}$}.
 			\end{cases} 
 		\end{align*}
 \end{thm}
 \begin{re}\label{re:rank-degree-bounds}
 	(i) By Theorem \ref{thm:class-number}, every norm Euclidean
 	$\mathcal{O}_{F}$-lattice is $\mathcal{O}_{F}$-maximal
 	and has class number one. Hanke and Kirschmer
     classified all positive definite
 	$\mathcal{O}_{F}$-maximal lattices of rank $\ge 3$ with
 	class number one over $F=\mathbb{Q}$ and over totally real fields $F\not=\mathbb{Q}$, respectively. See 
 	\cite[Theorem 7.5 and \S8]{hanke_maximal_2011} and 
 	\cite[Theorem 7.1]{kirschmer_one-class_2014}. By Theorem \ref{thm:class-number} and
 	Remark \ref{re:norm-bound-high-rank}, their lists provide all candidates for
 	positive definite norm Euclidean lattices of rank $\ge3$
 	over totally real number fields.
 	
 	(ii) By \cite[(1.2)]{he_trace_2026}, every positive definite
 	trace Euclidean lattice is norm Euclidean.
 	Consequently, Theorem \ref{thm:rank-degree-bounds} sharpens the
 	bounds in \cite[Corollary 4.3]{he_trace_2026}
 	for lattices of rank at least three.  
 \end{re}

The paper is organized as follows.
In Section \ref{sec:norm-descent}, we introduce some $F$-linear maps
arising from Eichler transformations and use them to establish
the descent of isotropic vectors. In Section \ref{sec:hyperbolic-cancellation}, we prove a cancellation
theorem for the hyperbolic lattice $\mathbf{H}$
and a related result involving $\mathbf{H}(\mathfrak{a})$.
In Section \ref{sec:class-number}, we relate lattice class numbers
to the class number of the base field and prove
Theorem \ref{thm:class-number} and
Corollary \ref{cor:stabilized-genus}. 
In Section \ref{sec:finiteness}, we prove
Theorems \ref{thm:finiteness} and \ref{thm:rank-degree-bounds}.

\section{$F$-linear maps and norm descent}
\label{sec:norm-descent}

Throughout this section, let $L$ be an integral
$\mathcal{O}_{F}$-lattice of rank $n$ on an anisotropic quadratic
space $V$ and put $M=L\perp \mathbf{H}$ as in \eqref{lorentzian}.

For an $\mathcal{O}_{F}$-lattice $J$ and $x\in FJ\setminus\{0\}$, 
recall from \cite[\S 81B]{omeara_quadratic_1963} that the coefficient of $x$ in $J$ is defined by  $\mathfrak{a}_{J}(x):=\{\alpha\in F:\alpha x\in J\}$, which is a nonzero fractional ideal. Then we write $	\mathfrak{c}_{J}(x)$ for the inverse of the \emph{coefficient} of $x$
in $J$:
\begin{align*}
	\mathfrak{c}_{J}(x)
	=\mathfrak{a}_{J}(x)^{-1}
	\quad\text{and}\quad \mathfrak{c}_{J}(0)=0.
\end{align*}
 We call $x\in J$ \textit{primitive} if
$\mathfrak{c}_{J}(x)=\mathcal{O}_{F}$.

\begin{prop}\label{prop:content-ideal}
	Let $J,J^{\prime}$ be $\mathcal{O}_{F}$-lattices, and let
	$x\in FJ$ and $x^{\prime}\in FJ^{\prime}$.
	\begin{enumerate}[label=\rm(\roman*), itemindent=0em, leftmargin=*]
		\item For every nonzero fractional ideal $\mathfrak{i}$, 
			$x\in\mathfrak{i}J$ if and only if
			$\mathfrak{c}_{J}(x)\subseteq\mathfrak{i}$. 
		\item For every $\alpha\in F$, $
			\mathfrak{c}_{J}(\alpha x)=\alpha\mathfrak{c}_{J}(x)$ and $\mathfrak{c}_{J\oplus J^{\prime}}(x,x^{\prime})
			=\mathfrak{c}_{J}(x)+\mathfrak{c}_{J^{\prime}}(x^{\prime}).$
		\item For every $\mathcal{O}_{F}$-module isomorphism
		$\sigma:J\to J^{\prime}$, extended $F$-linearly to $FJ$, we have
		$\mathfrak{c}_{J^{\prime}}(\sigma(x))=\mathfrak{c}_{J}(x)$.
	\end{enumerate}
\end{prop}

\begin{proof}
	By \cite[Theorem 81:3]{omeara_quadratic_1963}, we may write
	$J=\sum_{k=1}^{r}\mathfrak{a}_{k}v_{k}$,
	where $v_{1},\ldots,v_{r}$ form a basis of $FJ$ and
	the $\mathfrak{a}_{k}$ are nonzero fractional ideals.
	For $x=\sum_{k=1}^{r}x_{k}v_{k}$, we claim that
	\begin{align}\label{eq:content-coordinates}
		\mathfrak{c}_{J}(x)=\sum_{k=1}^{r}x_{k}\mathfrak{a}_{k}^{-1}.
	\end{align}
	If $x=0$, it is trivial. So we may assume that $x\neq0$. Then from
	\cite[Examples 22:4 and 81:4]{omeara_quadratic_1963}, we have
	\begin{align*}
		\mathfrak{c}_{J}(x)
		&=\left(\bigcap_{x_{k}\neq0}x_{k}^{-1}\mathfrak{a}_{k}\right)^{-1}
		=\sum_{k=1}^{r}x_{k}\mathfrak{a}_{k}^{-1}.
	\end{align*}
	Thus the claim is proved. 
	
	For (i), by \eqref{eq:content-coordinates}, we have the following equivalence:
	\begin{align*}
		x\in\mathfrak{i}J
		&\quad\Longleftrightarrow\quad
		x_{k}\in\mathfrak{i}\mathfrak{a}_{k}\text{ for every }k \quad\Longleftrightarrow\quad
		\sum_{k=1}^{r}x_{k}\mathfrak{a}_{k}^{-1}\subseteq\mathfrak{i}.
	\end{align*}
	For (ii), the identities follow from
	\eqref{eq:content-coordinates}.
	For (iii), the case $x=0$ follows from the definition.
	For $x\neq0$,  we see that
	\begin{align*}
		\mathfrak{a}_{J^{\prime}}(\sigma(x))
		&=\{\alpha\in F\mid\alpha\sigma(x)\in J^{\prime}\}=\{\alpha\in F\mid\sigma(\alpha x)\in\sigma(J)\}\\
		&=\{\alpha\in F\mid \alpha x\in J\}
		=\mathfrak{a}_{J}(x).
	\end{align*}
	Taking inverse ideals proves (iii).
\end{proof}
For $z=(v,\alpha,\beta)=v+\alpha e+\beta f\in M$, we have
\begin{align*}
	Q(z)=Q(v)+\alpha\beta
	\quad\text{and}\quad
	\mathfrak{c}_{M}(z)=
	\mathfrak{c}_{L}(v)
	+\alpha\mathcal{O}_{F}+\beta\mathcal{O}_{F},
\end{align*}
from the definition of $M$ and Proposition \ref{prop:content-ideal}(ii). 

As in
\cite[\S3]{omeara_generation_1965}, for a nonzero isotropic vector
$i\in FM$ and $w\in i^{\perp}\coloneq\{w\in FM\mid B(i,w)=0\}$,
we define the \textit{Eichler transformation} by
\begin{align}\label{eq:eichler}
	E_{w}^{i}(x)
	\coloneq x+B(x,i)w-B(x,w)i-\dfrac{1}{2}Q(w)B(x,i)i,
	\qquad x\in FM.
\end{align}
This is clearly an $F$-linear map. We summarize the identities from
\cite[\S3, (3.2)--(3.3)]{eichler_quadratische_1974} and include a proof for completeness.
\begin{prop}\label{prop:eichler-identities}
	Let $i\in FM\setminus\{0\}$ satisfy $Q(i)=0$.
	For $w,w_{1},w_{2}\in i^{\perp}$, we have
	 \begin{enumerate}[label=\rm(\roman*), itemindent=0em, leftmargin=*]
		\item $ E_{w}^{i}(i)=i$,  $E_{w_{1}}^{i}E_{w_{2}}^{i}=E_{w_{1}+w_{2}}^{i}$ and $(E_{w}^{i})^{-1}=E_{-w}^{i}$.
		
		\item   $E_{w}^{i}\in O(FM)$.
	\end{enumerate} 
\end{prop}
\begin{proof}
	(i) Since $B(i,i)=Q(i)=0$ and $B(i,w)=0$, we have
	\begin{align*}
		E_{w}^{i}(i)
		&=i+B(i,i)w-B(i,w)i-\dfrac{1}{2}Q(w)B(i,i)i
		=i.
	\end{align*}
	This proves the first equality.
	For $x\in FM$, put $
		\alpha=B(x,i)
		\quad\text{and}\quad
		\beta_{k}=B(x,w_{k})$, $k=1,2$.
	By \eqref{eq:eichler},
	\begin{align*}
		E_{w_{2}}^{i}(x)
		=x+\alpha w_{2}-\left(\beta_{2}+\dfrac{1}{2}\alpha Q(w_{2})\right)i.
	\end{align*}
	Since $Q(i)=B(i,w_{1})=B(i,w_{2})=0$, it follows that
	\begin{align*}
		B(E_{w_{2}}^{i}(x),i)
		&=\alpha+\alpha B(w_{2},i)
		-\left(\beta_{2}+\dfrac{1}{2}\alpha Q(w_{2})\right)Q(i)
		=\alpha,\\
		B(E_{w_{2}}^{i}(x),w_{1})
		&=\beta_{1}+\alpha B(w_{2},w_{1}) -\left(\beta_{2}+\dfrac{1}{2}\alpha Q(w_{2})\right)B(i,w_{1})=\beta_{1}+\alpha B(w_{2},w_{1}).
	\end{align*}
    One can check that
	\begin{align*}
		E_{w_{1}}^{i}E_{w_{2}}^{i}(x)
		&=E_{w_{2}}^{i}(x)+\alpha w_{1}-\left(\beta_{1}+\alpha B(w_{2},w_{1})
		+\dfrac{1}{2}\alpha Q(w_{1})\right)i\\
		&=x+\alpha(w_{1}+w_{2})-(\beta_{1}+\beta_{2})i-\dfrac{1}{2}\alpha
		\left(Q(w_{1})+2B(w_{1},w_{2})+Q(w_{2})\right)i.	
	\end{align*}
	Since  $w_{1}+w_{2}\in i^{\perp}$, $E_{w_{1}+w_{2}}^{i}(x)$ is defined, and we have
	\begin{align*}
			E_{w_{1}+w_{2}}^{i}(x)&=x+\alpha(w_{1}+w_{2})-\left(\beta_{1}+\beta_{2}
		+\dfrac{1}{2}\alpha Q(w_{1}+w_{2})\right)i.
	\end{align*}
	 So the second equality follows from the identity
	 \begin{align*}
	 	Q(w_{1}+w_{2})
	 	=Q(w_{1})+2B(w_{1},w_{2})+Q(w_{2}).
	 \end{align*}
	 Taking $(w_{1},w_{2})=(w,-w),(-w,w)$ in the second equality, we obtain
	\begin{align*}
		E_{w}^{i}E_{-w}^{i}
		=E_{-w}^{i}E_{w}^{i}=E_{0}^{i}=\id.
	\end{align*}
	Thus $(E_{w}^{i})^{-1}=E_{-w}^{i}$.

	(ii) For $x\in FM$, put $\alpha=B(x,i)$,
	$\beta=B(x,w)$ and $\gamma=\beta+\alpha Q(w)/2$.
	Then $E_{w}^{i}(x)=x+\alpha w-\gamma i$.
	Since $Q(i)=B(w,i)=0$, we have
	\begin{align*}
		Q(E_{w}^{i}(x))
		&=Q(x+\alpha w-\gamma i)\\
		&=Q(x)+\alpha^{2}Q(w)+\gamma^{2}Q(i)+2\alpha B(x,w)-2\gamma B(x,i)-2\alpha\gamma B(w,i)\\
		&=Q(x)+\alpha^{2}Q(w)+2\alpha\beta
		-2\alpha\left(\beta+\dfrac{1}{2}\alpha Q(w)\right)=Q(x).
	\end{align*}
	Combining this with the invertibility in (i), we deduce 
	$E_{w}^{i}\in O(FM)$.
\end{proof}
\begin{lem}\label{lem:elementary-isometries}
	For $u\in L$ and $\gamma\in\mathcal{O}_{F}^{\times}$, define the $F$-linear maps
	$T_{u},S,R_{\gamma}:FM\to FM$ by
	\begin{align}\label{eq:descent-maps}
		\begin{aligned}
			R_{\gamma}(v+\alpha e+\beta f)&\coloneq v+\gamma^{-1}\alpha e+\gamma\beta f,\\		
			S(v+\alpha e+\beta f)&\coloneq v+\beta e+\alpha f, \\
			T_{u}(v+\alpha e+\beta f)
			&\coloneq v-\alpha u+\alpha e+\bigl(\beta+2B(v,u)-\alpha Q(u)\bigr)f, 
		\end{aligned} 
	\end{align}
	where $v\in V$ and $\alpha,\beta\in F$. Then the following hold.
	\begin{enumerate}[label=\rm(\roman*), itemindent=0em, leftmargin=*]
		\item $R_{\gamma},S,T_{u}\in O(M)$.
		\item $R_{\gamma}(\gamma e)=e$, $ST_{u}S=E_{-2u}^{e}$,
			$ST_{u}S(e)=e$ and $ST_{u}S(u-Q(u)e+f)=f$.
		\item $\mathfrak{c}_{M}(T_{u}(z))=\mathfrak{c}_{M}(z)$ for any $z\in M$.
	\end{enumerate}
\end{lem}

\begin{proof}
	(i) Since $\gamma,\gamma^{-1}\in\mathcal{O}_{F}^{\times}$,
	$R_{\gamma}(M)=M$. By \eqref{eq:descent-maps},
	\begin{align*}
		R_{\gamma^{-1}}R_{\gamma}(v+\alpha e+\beta f)=R_{\gamma^{-1}}(v+\gamma^{-1}\alpha e+\gamma\beta f)=v+\alpha e+\beta f.
	\end{align*}
	Replacing $\gamma$ by $\gamma^{-1}$ gives
	$R_{\gamma}R_{\gamma^{-1}}=\id$ as well.
	Thus $R_{\gamma}^{-1}=R_{\gamma^{-1}}$.
	Moreover, for $v\in V$ and $\alpha,\beta\in F$,
	\begin{align*}
		Q(R_{\gamma}(v+\alpha e+\beta f))
		&=Q(v)+(\gamma^{-1}\alpha)(\gamma\beta)=Q(v)+\alpha\beta=Q(v+\alpha e+\beta f).
	\end{align*}
	Hence $R_{\gamma}\in O(M)$.

	By definition, we have $
		S(M)=M$ and $S^{-1}=S$.
	Moreover,
	\begin{align*}
		Q(S(v+\alpha e+\beta f))=Q(v+\beta e+\alpha f)=Q(v)+\beta\alpha=Q(v+\alpha e+\beta f).
	\end{align*}
	Therefore, $S\in O(M)$.

	Let $x=v+\alpha e+\beta f\in FM$. Since $V\perp(Fe+Ff)$ and
	$B(e,f)=1/2$, we have
	\begin{align*}
		B(x,f)&=\dfrac{\alpha}{2},\quad B(x,-2u)=-2B(v,u),\quad Q(-2u)=4Q(u).
	\end{align*}
	Since $Q(f)=B(f,u)=0$, the transformation $E_{-2u}^{f}$
	is defined, and \eqref{eq:eichler} gives
	\begin{align*}
		E_{-2u}^{f}(x)&=x+B(x,f)(-2u)-B(x,-2u)f-\dfrac{1}{2}Q(-2u)B(x,f)f\\
		&=x-\alpha u+2B(v,u)f
		-\alpha Q(u)f\\
		&=v-\alpha u+\alpha e+\bigl(\beta+2B(v,u)-\alpha Q(u)\bigr)f
		=T_{u}(x).
	\end{align*}
	By Proposition \ref{prop:eichler-identities}(i) and (ii), $
		T_{u}=E_{-2u}^{f}\in O(FM)$ and $T_{u}^{-1}=E_{2u}^{f}=T_{-u}$. 
	Let $v\in L$ and $\alpha,\beta\in\mathcal{O}_{F}$.
	Then $v-\alpha u\in L$. Since $L$ is integral,
	\begin{align*}
		Q(u)\in\mathcal{O}_{F}
		\quad\text{and}\quad
		2B(v,u)=Q(v+u)-Q(v)-Q(u)\in\mathcal{O}_{F}.
	\end{align*}
	Hence $\beta+2B(v,u)-\alpha Q(u)\in\mathcal{O}_{F}$.
	It follows from \eqref{eq:descent-maps} that $T_{u}(M)\subseteq M$.
	Applying the same argument to $-u\in L$, we have
	$T_{u}^{-1}(M)=T_{-u}(M)\subseteq M$, whence $M\subseteq T_{u}(M)$.
	Therefore, $T_{u}(M)=M$ and $T_{u}\in O(M)$. (i) is proved.

	(ii) By definition, we have $R_{\gamma}(\gamma e)=\gamma R_{\gamma}(e)=e$.
	Let $x=v+\alpha e+\beta f\in FM$, and write
	$x=(v,\alpha,\beta)$. We obtain
	\begin{align*}
		ST_{u}S(v,\alpha,\beta)
		=ST_{u}(v,\beta,\alpha)&=S\bigl(v-\beta u,\ \beta,\ \alpha+2B(v,u)-\beta Q(u)\bigr)\\
		&=\bigl(v-\beta u,\ \alpha+2B(v,u)-\beta Q(u),\ \beta\bigr).
	\end{align*}
	Since $B(x,e)=\beta/2$, we have
	\begin{align*}
		E_{-2u}^{e}(x)
		&=x-\beta u+\bigl(2B(v,u)-\beta Q(u)\bigr)e=\bigl(v-\beta u,\ \alpha+2B(v,u)-\beta Q(u),\ \beta\bigr).
	\end{align*}
	Hence $ST_{u}S=E_{-2u}^{e}$. By Proposition
	\ref{prop:eichler-identities}(i), $ST_{u}S(e)=E_{-2u}^{e}(e)=e$.
	Since $B(u,u)=Q(u)$, by \eqref{eq:descent-maps}, one can verify that
	\begin{align*}
		ST_{u}S(u-Q(u)e+f)
		&=ST_{u}(u+e-Q(u)f)\\
		&=S\bigl(u-u+e+\bigl(-Q(u)+2B(u,u)-Q(u)\bigr)f\bigr)=S(e)=f.
	\end{align*}

	(iii) By (i), $T_{u}$ is an $\mathcal{O}_{F}$-module
	automorphism of $M$, so the assertion follows from
	Proposition \ref{prop:content-ideal}(iii).
\end{proof}

\begin{lem}\label{lem:descent}
	Suppose that $L$ is norm Euclidean. For every
	$z\in M\setminus\{0\}$ with $Q(z)=0$, there exist
	$\sigma\in O(M)$ and $\gamma\in\mathcal{O}_{F}\setminus\{0\}$ such that
	\begin{align*}
		\sigma(z)=\gamma e
		\quad\text{and}\quad
		\mathfrak{c}_{M}(z)=\gamma\mathcal{O}_{F}.
	\end{align*}
	Moreover,
	\begin{align*}
		\{z\in M\mid Q(z)=0\;\mbox{and}\; \mathfrak{c}_{M}(z)=\mathcal{O}_{F}\}
		=O(M)e.
	\end{align*}
\end{lem}

\begin{proof}
	Let $z\in M\setminus\{0\}$ satisfy $Q(z)=0$, and write
	$z=(w,\alpha_{0},\beta_{0})$ with $w\in L$ and
	$\alpha_{0},\beta_{0}\in\mathcal{O}_{F}$.
	If $\alpha_{0}=\beta_{0}=0$, then $Q(w)=Q(z)=0$. Since $V$ is
	anisotropic, $w=0$, contrary to $z\neq0$.
	Thus at least one of $\alpha_{0}$ and $\beta_{0}$ is nonzero.
	If $\alpha_{0}=0$, then by Lemma \ref{lem:elementary-isometries}(i),
	\begin{align*}
		S(z)=(w,\beta_{0},0)\in O(M)z
		\quad\text{and}\quad \beta_{0}\neq0.
	\end{align*}
	Hence $O(M)z$ contains a vector with nonzero $e$-coordinate.
	Every vector in this orbit lies in $M$, so its two hyperbolic
	coordinates belong to $\mathcal{O}_{F}$. Set
	\begin{align*}
		\mathcal{N}
		\coloneq\{|\nr_{F/\mathbb{Q}}(\alpha)|:
		(v,\alpha,\beta)\in O(M)z,\ \alpha\neq0\}.
	\end{align*}
	Since $\emptyset\neq\mathcal{N}\subseteq\mathbb{Z}_{>0}$, the set
	$\mathcal{N}$ has a minimum. Choose $\sigma_{0}\in O(M)$ such that
	$\sigma_{0}(z)=(v,\alpha,\beta)$ with $\alpha\neq 0$ and
	$|\nr_{F/\mathbb{Q}}(\alpha)|=\min\mathcal{N}$. Write $z_{0}=\sigma_{0}(z)$.
	Since $\sigma_{0}$ is an isometry,
	\begin{align*}
		Q(v)+\alpha\beta=Q(z_{0})=Q(z)=0,
	\end{align*}
	and hence $\beta=-Q(v)/\alpha$. For any $u\in L$, by \eqref{eq:descent-maps} and
	Lemma \ref{lem:elementary-isometries}(i), we have
	$T_{u}(z_{0})=(v^{\prime},\alpha,\beta^{\prime})\in M$ with 
	\begin{align*}
		v^{\prime}=v-\alpha u \in L
		\quad\text{and}\quad
		\beta^{\prime}=\beta+2B(v,u)-\alpha Q(u)\in \mathcal{O}_{F}.
	\end{align*}
	Substituting $\beta=-Q(v)/\alpha$, we obtain
	\begin{align*}
		\beta^{\prime}
		&=-\dfrac{Q(v)}{\alpha}+2B(v,u)-\alpha Q(u)=-\alpha Q\left(\dfrac{v}{\alpha}-u\right).
	\end{align*}
	Since $v/\alpha\in V$ and $L$ is norm Euclidean, we may choose
	$u\in L$ such that
	\begin{align*}
		\left|\nr_{F/\mathbb{Q}}
		\left(Q\left(\dfrac{v}{\alpha}-u\right)\right)\right|<1.
	\end{align*}
	By multiplicativity of $\nr_{F/\mathbb{Q}}$,
	\begin{align*}
		|\nr_{F/\mathbb{Q}}(\beta^{\prime})|
		&=|\nr_{F/\mathbb{Q}}(\alpha)|
		\left|\nr_{F/\mathbb{Q}}
		\left(Q\left(\dfrac{v}{\alpha}-u\right)\right)\right|<|\nr_{F/\mathbb{Q}}(\alpha)|.
	\end{align*}
	By Lemma \ref{lem:elementary-isometries}(i), we have
	\begin{align*}
		ST_{u}(z_{0})=(v^{\prime},\beta^{\prime},\alpha)\in O(M)z.
	\end{align*}
	If $\beta^{\prime}\neq0$, then
	$|\nr_{F/\mathbb{Q}}(\beta^{\prime})|\in\mathcal{N}$, which contradicts the choice of $\alpha$. So $\beta^{\prime}=0$.
	Since $\beta^{\prime}=0$ and $T_{u}\in O(M)$,
	\begin{align*}
		0=Q(z_{0})=Q(T_{u}(z_{0}))
		=Q(v^{\prime})+\alpha\beta^{\prime}=Q(v^{\prime}).
	\end{align*}
	Since $V$ is anisotropic, $v^{\prime}=0$ and thus $ T_{u}(z_{0})=(0,\alpha,0)=\alpha e$. Take $\sigma=T_{u}\sigma_{0}\in O(M)$ and
	$\gamma=\alpha\in\mathcal{O}_{F}\setminus\{0\}$.
	Then $\sigma(z)=T_{u}\sigma_{0}(z)=T_{u}(z_{0})=\gamma e$, as required.

	By Proposition \ref{prop:content-ideal}(iii) and
	Lemma \ref{lem:elementary-isometries}(iii), we have
	\begin{align*}
		\mathfrak{c}_{M}(z)
		=\mathfrak{c}_{M}(\sigma_{0}(z))	=\mathfrak{c}_{M}(z_{0})
		=\mathfrak{c}_{M}(T_{u}(z_{0}))
		=\mathfrak{c}_{M}(\gamma e).
	\end{align*}
Also, by Proposition \ref{prop:content-ideal}(ii), $
		\mathfrak{c}_{M}(\gamma e)
		=\mathfrak{c}_{L}(0)+\gamma\mathcal{O}_{F}
		+0\mathcal{O}_{F}=\gamma\mathcal{O}_{F}$. Hence $\mathfrak{c}_{M}(z)=\gamma\mathcal{O}_{F}$.
	This proves the first assertion.

	For the last assertion, let $z\in M$ satisfy $Q(z)=0$ and
	$\mathfrak{c}_{M}(z)=\mathcal{O}_{F}$. Since $\mathfrak{c}_{M}(0)=0$, we have
	$z\neq0$, so the first assertion applies. With $\sigma$ and $\gamma$
	as above, we have $\gamma\mathcal{O}_{F}=\mathfrak{c}_{M}(z)=\mathcal{O}_{F}$ and hence
	$\gamma\in\mathcal{O}_{F}^{\times}$.
	By Lemma \ref{lem:elementary-isometries}(i) and (ii), we have
	$R_{\gamma}\in O(M)$. Hence
	\begin{align*}
		R_{\gamma}\sigma(z)=R_{\gamma}(\gamma e)=e
		\quad\text{and}\quad
		z=(R_{\gamma}\sigma)^{-1}(e)\in O(M)e.
	\end{align*}
	
	Conversely, let $z=\tau(e)$ with $\tau\in O(M)$. Then $	Q(z)=Q(e)=0$, and
	Proposition \ref{prop:content-ideal}(iii) implies that $
		\mathfrak{c}_{M}(z)=\mathfrak{c}_{M}(e)=\mathcal{O}_{F}$.
\end{proof} 

\begin{lem}\label{lem:isotropic-prime}
	Suppose that $L$ is norm Euclidean, and let $\mathfrak{p}$ be a nonzero prime ideal of $\mathcal{O}_{F}$.
	If the quadratic form $\overline{Q}$ induced by $Q$ on $L/\mathfrak{p}L$ is isotropic, then $\mathfrak{p}$ is principal.
\end{lem}

\begin{proof}
	Let $\overline{v}$ be
	a nonzero isotropic vector of the induced form $\overline{Q}$ on $L/\mathfrak{p}L$. Choose $\alpha\in\mathfrak{p}\setminus\mathfrak{p}^{2}$ and put
	$\mathfrak{i}=\alpha\mathfrak{p}^{-1}$.
	Then $\mathfrak{i}$ is a nonzero integral ideal satisfying
	\begin{align*}
		\mathfrak{i}\mathfrak{p}=\alpha\mathcal{O}_{F}
		\quad\text{and}\quad \mathfrak{i}+\mathfrak{p}=\mathcal{O}_{F}.
	\end{align*}
	Since $\mathfrak{i}L+\mathfrak{p}L=L$, reduction maps
	$\mathfrak{i}L$ onto $L/\mathfrak{p}L$.
	Choose $v\in\mathfrak{i}L$ lifting $\overline{v}$.
	Then $v\notin\mathfrak{p}L$.
	
	Since $L$ is integral and $v\in \mathfrak{i}L$, $
	Q(v)\in\mathfrak{i}^{2}$.  Since $\overline{Q}(\overline{v})=0$, $Q(v)\in\mathfrak{p}$.
	These relations, combined with $\mathfrak{i}^{2}+\mathfrak{p}=\mathcal{O}_{F}$, show that
	\begin{align*}
		Q(v)\in\mathfrak{i}^{2}\cap\mathfrak{p}
		=\mathfrak{i}^{2}\mathfrak{p}
		=\alpha\mathfrak{i}.
	\end{align*}
	Put $\beta=-Q(v)/\alpha\in\mathfrak{i}$ and
	$z=(v,\alpha,\beta)\in M$. Since
	$v\in\mathfrak{i}L$, Proposition \ref{prop:content-ideal}(i) implies
	$\mathfrak{c}_{L}(v)\subseteq\mathfrak{i}$.
	Applying Proposition \ref{prop:content-ideal}(ii), we have
	\begin{align*}
		\mathfrak{c}_{M}(z)
		=\mathfrak{c}_{L}(v)+\alpha\mathcal{O}_{F}
		+\beta\mathcal{O}_{F}.
	\end{align*}
	Since $\alpha,\beta\in\mathfrak{i}$, it follows that $
	\alpha\mathcal{O}_{F}\subseteq\mathfrak{c}_{M}(z)
	\subseteq\mathfrak{i}$.
	So, for every prime ideal $\mathfrak{q}$,
	\begin{align*}
		\ord_{\mathfrak{q}}(\mathfrak{i})
		\leq\ord_{\mathfrak{q}}(\mathfrak{c}_{M}(z))
		\leq\ord_{\mathfrak{q}}(\alpha).
	\end{align*}
	If $\mathfrak{q}\neq\mathfrak{p}$, since
	$\alpha\mathcal{O}_{F}=\mathfrak{i}\mathfrak{p}$,
	$\ord_{\mathfrak{q}}(\alpha)=\ord_{\mathfrak{q}}(\mathfrak{c}_{M}(z))
	=\ord_{\mathfrak{q}}(\mathfrak{i})$. By Proposition \ref{prop:content-ideal}(i) and
	$v\notin\mathfrak{p}L$, we have
	$\mathfrak{c}_{L}(v)\not\subseteq\mathfrak{p}$.
	Since $\mathfrak{c}_{L}(v)\subseteq\mathfrak{c}_{M}(z)$, we obtain
	\begin{align*}
		\ord_{\mathfrak{p}}(\mathfrak{c}_{M}(z))=0
		=\ord_{\mathfrak{p}}(\mathfrak{i}),
	\end{align*}
	where the last equality follows from $\mathfrak{p}\nmid\mathfrak{i}$.
	Thus $\mathfrak{c}_{M}(z)$ and $\mathfrak{i}$ have the same valuation
	at every finite prime, and hence $\mathfrak{c}_{M}(z)=\mathfrak{i}$.
	
	By construction, $z\neq0$ and $Q(z)=Q(v)+\alpha\beta=0$.
	By Lemma \ref{lem:descent}, there exists
	$\gamma\in\mathcal{O}_{F}\setminus\{0\}$ such that
	\begin{align*}
		\mathfrak{i}=\mathfrak{c}_{M}(z)=\gamma\mathcal{O}_{F}.
	\end{align*}
	Hence $\mathfrak{p}=(\alpha/\gamma)\mathcal{O}_{F}$ is principal.
\end{proof}

\section{Cancellation of $\mathbf{H}$}
\label{sec:hyperbolic-cancellation}

Throughout this section, let $F$ be an algebraic number field.
Let $\Omega_{F}$ be the set of all primes of $F$, and let
$\infty_{F}$ be the set of archimedean primes. Let $V$ be a non-degenerate anisotropic quadratic space over $F$ and $L$ be an 
$\mathcal{O}_{F}$-lattice of rank $n$ on $V$. Put $M=L\perp\mathbf{H}$, and retain the notation
$e,f$ and $ \mathfrak{c}_{M}(z)$ from
Section \ref{sec:norm-descent}.

 \begin{thm}\label{thm:cancellation}
 	Suppose that $L$ is norm Euclidean.
 	If $L^{\prime}$ is an $\mathcal{O}_{F}$-lattice and
 	$L^{\prime}\perp\mathbf{H}\cong L\perp\mathbf{H}$,
 	then $L^{\prime}\cong L$.
 \end{thm}
 
 \begin{proof}
 	Write $M^{\prime}=L^{\prime}\perp\mathbf{H}^{\prime}$, where
 	$\mathbf{H}^{\prime}=\mathcal{O}_{F}e^{\prime}
 	+\mathcal{O}_{F}f^{\prime}$ is the hyperbolic summand
 	with its standard basis. Let $\varphi:M^{\prime}\to M$ be an isometry. Note that for any $\gamma\in F$, $\gamma e^{\prime}\in M^{\prime}$ if and only if $\gamma\in \mathcal{O}_{F}$. So $\mathfrak{a}_{M^{\prime}}(e^{\prime})=\mathcal{O}_{F}$ and hence
 	\begin{align*}
 		\mathfrak{c}_{M^{\prime}}(e^{\prime})=\mathfrak{a}_{M^{\prime}}(e^{\prime})^{-1}=\mathcal{O}_{F}.
 	\end{align*}
 	This combined with Proposition \ref{prop:content-ideal}(iii) shows that
 	$\mathfrak{c}_{M}(\varphi(e^{\prime}))
 	=\mathfrak{c}_{M^{\prime}}(e^{\prime})=\mathcal{O}_{F}$.
 	Since $\varphi$ is an isometry, we also have
 	$Q(\varphi(e^{\prime}))=Q(e^{\prime})=0$. Therefore, $\varphi(e^{\prime})$ is primitive and isotropic.
 	Since $L$ is norm Euclidean, by the last assertion of
 	Lemma \ref{lem:descent}, there exists $\sigma\in O(M)$ such that
 	$\sigma\varphi(e^{\prime})=e$. Put $\psi=\sigma\varphi$.
 	
 	Write $\psi(f^{\prime})=w+\alpha e+\beta f$, where
 	$w\in L$ and $\alpha,\beta\in\mathcal{O}_{F}$.
 	Since $\psi$ is an isometry and $\psi(e^{\prime})=e$, we have
 	\begin{align*}
 		\dfrac{\beta}{2}
 		=B(e,\psi(f^{\prime}))
 		=B(\psi(e^{\prime}),\psi(f^{\prime}))
 		=B(e^{\prime},f^{\prime})=\dfrac{1}{2}.
 	\end{align*}
 	Hence $\beta=1$. Also,
 	\begin{align*}
 		0=Q(f^{\prime})=Q(\psi(f^{\prime}))
 		=Q(w)+\alpha\beta=Q(w)+\alpha.
 	\end{align*}
 	Thus, $\alpha=-Q(w)$ and so
 	\begin{align*}
 		\psi(f^{\prime})=w-Q(w)e+f.
 	\end{align*}
 	By Lemma \ref{lem:elementary-isometries}(i) and (ii), we have $E_{-2w}^{e}=ST_{w}S\in O(M)$,
 	\begin{align*} 
 			E_{-2w}^{e}\psi(e^{\prime})=E_{-2w}^{e}(e)=e \quad\text{and}\quad E_{-2w}^{e}\psi(f^{\prime})=ST_{w}S(w-Q(w)e+f)=f.
 	\end{align*}
 	Therefore, $\Psi=E_{-2w}^{e}\psi:M^{\prime}\to M$ is an isometry satisfying
 	\begin{align*}
 		\Psi(e^{\prime})=e
 		\quad\text{and}\quad \Psi(f^{\prime})=f.
 	\end{align*}
 	
 	For $x=v+\alpha e+\beta f\in M$, we have
 	\begin{align*}
 		B(x,e)=\dfrac{\beta}{2}
 		\quad\text{and}\quad B(x,f)=\dfrac{\alpha}{2}.
 	\end{align*}
 	Hence $x\in L$ if and only if $B(x,e)=B(x,f)=0$.
 	The same argument for $M^{\prime}$ shows that
 	$x^{\prime}\in L^{\prime}$ if and only if
 	$B(x^{\prime},e^{\prime})=B(x^{\prime},f^{\prime})=0$.
 	Since the isometry $\Psi$ preserves $B$, for every $x^{\prime}\in M^{\prime}$,
 	\begin{align*}
 		x^{\prime}\in L^{\prime}
 		&\quad\Longleftrightarrow\quad
 		B(x^{\prime},e^{\prime})=B(x^{\prime},f^{\prime})=0\\
 		&\quad\Longleftrightarrow\quad
 		B(\Psi(x^{\prime}),e)=B(\Psi(x^{\prime}),f)=0 \quad\Longleftrightarrow\quad \Psi(x^{\prime})\in L.
 	\end{align*}
 	Since $\Psi$ maps $\mathbf{H}^{\prime}$ onto $\mathbf{H}$,
 	\begin{align*}
 		\Psi(L^{\prime})
 		&=\Psi\bigl(M^{\prime}\cap
 		(F\mathbf{H}^{\prime})^{\perp}\bigr)=M\cap(F\mathbf{H})^{\perp}=L,
 	\end{align*}
 	as desired.
 \end{proof}

\begin{cor}\label{cor:ideal-cancellation}
	 Let $J$ be an
	$\mathcal{O}_{F}$-lattice on a non-degenerate quadratic space and let
	$\mathfrak{a}$ be a nonzero fractional ideal of $F$. Suppose that $L$ is norm Euclidean. Then
	\begin{align*}
		J\perp\mathbf{H}(\mathfrak{a})\cong L\perp\mathbf{H}
		\quad\Longleftrightarrow\quad
		[\mathfrak{a}]=1\ \text{in }\cl(F)
		\quad\text{and}\quad J\cong L.
	\end{align*}
\end{cor}

\begin{proof}
	Sufficiency is clear. For necessity, put $N=J\perp\mathbf{H}(\mathfrak{a})$ and suppose that
	$\varphi:N\to M=L\perp\mathbf{H}$ is an isometry.
	Choose $0\neq \alpha\in\mathfrak{a}$ and put $z=\varphi(\alpha e)$.
	Then $z\neq0$, $Q(z)=0$, and
	\begin{align*}
		\mathfrak{a}_{N}(\alpha e)
		&=\{\lambda\in F:\lambda\alpha\in\mathfrak{a}\}
		=\alpha^{-1}\mathfrak{a}.
	\end{align*}
	So $\mathfrak{c}_{N}(\alpha e)=\alpha\mathfrak{a}^{-1}$.         By Proposition \ref{prop:content-ideal}(iii) and
	Lemma \ref{lem:descent}, there is
	$\gamma\in\mathcal{O}_{F}\setminus\{0\}$ such that
	\begin{align*}
		\alpha\mathfrak{a}^{-1}
		=\mathfrak{c}_{N}(\alpha e)=\mathfrak{c}_{M}(\varphi(\alpha e))=\mathfrak{c}_{M}(z)
		=\gamma\mathcal{O}_{F}.
	\end{align*}
	Hence $\mathfrak{a}=(\alpha/\gamma)\mathcal{O}_{F}$ is principal.
	Write $\mathfrak{a}=\beta\mathcal{O}_{F}$ with
	$\beta\in F^{\times}$. Then $\mathbf{H}(\mathfrak{a})=\mathcal{O}_{F}(\beta e)+\mathcal{O}_{F}(\beta^{-1}f)$. So the map given by
	\begin{align*}
		\beta e\mapsto e
		\quad\text{and}\quad \beta^{-1}f\mapsto f
	\end{align*}
	is an isometry $\mathbf{H}(\mathfrak{a})\cong\mathbf{H}$.
	Thus $J\perp\mathbf{H}\cong L\perp\mathbf{H}$, and
	Theorem \ref{thm:cancellation} implies $J\cong L$.
\end{proof}

 \begin{re}
 	When $F=\mathbb{Q}$, Clark and Jagy proved that nondegenerate
 	integral $\mathbb{Z}$-lattices $J$ and $J^{\prime}$ on the same
 	quadratic space over $\mathbb{Q}$ satisfy
 	\begin{align*}
 		J^{\prime}\in\gen(J)
 		\quad\Longleftrightarrow\quad
 		J^{\prime}\perp\mathbf{H}\cong J\perp\mathbf{H};
 	\end{align*}
 	see \cite[Theorem 6.3]{clark_ADC-II-2014}.
 	If $J$ is norm Euclidean, then by Theorem
 	\ref{thm:cancellation}, we have the further implication
 	\begin{align*}
 		J^{\prime}\in\gen(J)
 		\quad\Longrightarrow\quad
 		J^{\prime}\perp\mathbf{H}\cong J\perp\mathbf{H}
 		\quad\Longrightarrow\quad J^{\prime}\cong J.
 	\end{align*}
 \end{re}

\section{Class numbers of norm Euclidean lattices}
\label{sec:class-number}

In this section, we retain the notation and assumptions of
Section \ref{sec:hyperbolic-cancellation}. Also, for $\mathfrak{p}\in\Omega_{F}$, write $F_{\mathfrak{p}}$ for
the completion of $F$ and put
$V_{\mathfrak{p}}=V\otimes_{F}F_{\mathfrak{p}} $.
For $\mathfrak{p}\in \Omega_{F}\setminus \infty_{F}$, let
$\mathcal{O}_{F_{\mathfrak{p}}}$ be the valuation ring of
$F_{\mathfrak{p}}$, and put
\begin{align*}
	\kappa_{\mathfrak{p}}=\mathcal{O}_{F}/\mathfrak{p}
	\quad\text{and}\quad
	L_{\mathfrak{p}}=L\otimes_{\mathcal{O}_{F}}\mathcal{O}_{F_{\mathfrak{p}}}.
\end{align*}
Write $\theta$ for the spinor norm
\cite[\S55]{omeara_quadratic_1963}.
Let $\mathbb{I}_{F}$ be the idele group of $F$, and put
\begin{align*}
	\mathbb{U}_{F}=\left(\prod_{\mathfrak{p}\in\infty_{F}}F_{\mathfrak{p}}^{\times}\right)
	\times\prod_{\mathfrak{p}\nmid\infty}\mathcal{O}_{F_{\mathfrak{p}}}^{\times}.
\end{align*}
Denote by $\mathbb{I}_{M}$ the subgroup of ideles whose $\mathfrak{p}$-component
lies in $\theta(O^{+}(M_{\mathfrak{p}}))$ for every finite prime $\mathfrak{p}$.

 We relate $h^{+}(M)$ to $h(F)$ by slightly modifying
 \cite[Proposition 3.2]{xu-zhang_indefinite_2022}.

\begin{lem}\label{lem:hyperbolic-spinor}
	Suppose that $L$ is integral and $n\ge1$. Put $M=L\perp \mathbf{H} $. Then 
	\begin{enumerate}[label=\rm(\roman*), itemindent=0em, leftmargin=*]
		\item  For $\mathfrak{p}\in \Omega_{F}\setminus \infty_{F}$, 
		$\theta(O^{+}(M_{\mathfrak{p}}))
		\supseteq\mathcal{O}_{F_{\mathfrak{p}}}^{\times}F_{\mathfrak{p}}^{\times2}$; thus $\mathbb{U}_{F}\mathbb{I}_{F}^{2}\subseteq\mathbb{I}_{M}$. 
		
		\item For every $N\in \gen(M)$, $\spn^{+}(N)=\cls^{+}(N)$.
		
		\item  $h^{+}(M)\mid h(F)$.
	\end{enumerate} 
\end{lem} 

\begin{proof}
	(i) The hyperbolic plane $\mathbf{H}_{\mathfrak{p}}$ is universal at every finite
	prime $\mathfrak{p}$ by
	\cite[Proposition 2.3 and Corollary 2.9]{xu-zhang_indefinite_2022},
	and at every infinite prime by
	\cite[Theorem 2.3]{he-hu-xu_indefinite_2023} with $k=1$. So for every finite prime $\mathfrak{p}$, we have
	$Q(\mathbf{H}_{\mathfrak{p}})=\mathcal{O}_{F_{\mathfrak{p}}}$.
	Since $M_{\mathfrak{p}}$ splits off
	$\mathbf{H}_{\mathfrak{p}}$, we have
	\begin{align*}
		\theta(O^{+}(M_{\mathfrak{p}}))\supseteq 	\theta(O^{+}(\mathbf{H}_{\mathfrak{p}}))
		\supseteq\mathcal{O}_{F_{\mathfrak{p}}}^{\times}F_{\mathfrak{p}}^{\times2},
	\end{align*}
	and $\mathbb{U}_{F}\mathbb{I}_{F}^{2}\subseteq\mathbb{I}_{M}$ follows.
	
	(ii) Since $\rank\,M=n+2\ge3$ and $FM$ is isotropic at every $\mathfrak{p}\in \infty_{F}$, the assertion follows from \cite[104:5]{omeara_quadratic_1963}.
	
	(iii) Since $FM$ is isotropic,	by \cite[55:2a]{omeara_quadratic_1963},
	\begin{align*}
		\theta(O^{+}(FM))=F^{\times}/F^{\times2}.
	\end{align*}
	Hence, by (ii) and \cite[102:7]{omeara_quadratic_1963},
	\begin{align}\label{h+Mdivide}
		h^{+}(M)=g^{+}(M)
		=[\mathbb{I}_{F}:F^{\times}\mathbb{I}_{M}].
	\end{align}
	
	On the other hand, by \cite[33:14]{omeara_quadratic_1963}, $\mathbb{I}_{F}/F^{\times}\mathbb{U}_{F}\cong \cl(F)$ and so 
	\begin{align}\label{clf2}
		\mathbb{I}_{F}/ F^{\times}\mathbb{U}_{F}\mathbb{I}_{F}^{2}\cong \cl(F)/\cl(F)^{2}.
	\end{align}
	By (i), $\mathbb{U}_{F}\mathbb{I}_{F}^{2}\subseteq\mathbb{I}_{M}$. 
	This combined with \eqref{clf2} shows that
	\begin{align}\label{hFdivide}   
		[\mathbb{I}_{F}:F^{\times}\mathbb{I}_{M}]\mid 	[\mathbb{I}_{F}:F^{\times}\mathbb{U}_{F}\mathbb{I}_{F}^{2}] =[\cl(F):\cl(F)^{2}]\mid h(F).
	\end{align}
	Combining \eqref{h+Mdivide} and \eqref{hFdivide}, we deduce that $h^{+}(M)\mid h(F)$.
\end{proof}

\begin{lem}\label{lem:stabilized-class-number}
	Suppose that $L$ is norm Euclidean of rank $n\ge1$
	and $h(F)\leq2$. Then
	\begin{align*}
		h(L\perp\mathbf{H})=h(F).
	\end{align*}
\end{lem}
\begin{proof}
	Put $M=L\perp\mathbf{H}$.
	By assumption, $h(F)\in\{1,2\}$. Recall from \cite[102:2 Example]{omeara_quadratic_1963} that every class contains either one or two proper classes. So
	Lemma \ref{lem:hyperbolic-spinor}(iii) implies that
	\begin{align*}
		1\le h(M)\le h^{+}(M)\le h(F).
	\end{align*}
	Thus $h(M)=h(F)$ if $h(F)=1$.
	
	Suppose that $h(F)=2$. Let $\mathfrak{a}$ be any nonprincipal
	fractional ideal of $F$ and put
	$M_{\mathfrak{a}}=L\perp\mathbf{H}(\mathfrak{a})$.
	For every finite prime $\mathfrak{p}$, write
	$\mathfrak{a}_{\mathfrak{p}}
	=\beta_{\mathfrak{p}}\mathcal{O}_{F_{\mathfrak{p}}}$ with
	$\beta_{\mathfrak{p}}\in F_{\mathfrak{p}}^{\times}$.
	The map given by
	\begin{align*}
		e\longmapsto\beta_{\mathfrak{p}}e
		\quad\text{and}\quad
		f\longmapsto\beta_{\mathfrak{p}}^{-1}f
	\end{align*}
	is an isometry of determinant $1$ from $\mathbf{H}_{\mathfrak{p}}$
	onto $\mathbf{H}(\mathfrak{a})_{\mathfrak{p}}$.
	Extending it by the identity on $V_{\mathfrak{p}}$, we see that
	$(M_{\mathfrak{a}})_{\mathfrak{p}}\cong M_{\mathfrak{p}}$.
	Since $M_{\mathfrak{a}}$ and $M$ lie on $V\perp\mathbb{H}$,
	we have $M_{\mathfrak{a}}\in\gen(M)$.
	By Corollary \ref{cor:ideal-cancellation},
	$M_{\mathfrak{a}}\not\cong M$, since $\mathfrak{a}$ is nonprincipal. Hence
	\begin{align*}
		2\le h(M)\le h^{+}(M)\le h(F)=2.
	\end{align*}
	Thus $h(F)=h(M)$.
\end{proof}

\begin{lem}\label{lem:binary-stabilization}
	Let $L$ be an integral $\mathcal{O}_{F}$-lattice on an
	anisotropic binary quadratic space $V$, and put
	$E=F(\sqrt{-\det V})$.
	If $E=H_{F}$, then every $L^{\prime}\in\gen(L)$ satisfies
	\begin{align*}
		L^{\prime}\perp\mathbf{H}\cong L\perp\mathbf{H}.
	\end{align*}
\end{lem}
\begin{proof}
	Let $L^{\prime}\in\gen(L)$ and put
	$M=L\perp\mathbf{H}$ and
	$M^{\prime}=L^{\prime}\perp\mathbf{H}$.
	Since $V$ is anisotropic and $E=H_{F}$, by
	\cite[Theorem 5.23]{cox_primes_2013}, we have
	\begin{align*}
		h(F)=[H_{F}:F]=[E:F]=2.
	\end{align*}
	The set
	\begin{align*}
		\mathcal{T}\coloneq\{\mathfrak{p}\in\Omega_{F}\setminus\infty_{F}:
		L_{\mathfrak{p}}\neq L^{\prime}_{\mathfrak{p}}\}
	\end{align*}
	is finite. By \cite[91:4a and \S102A]{omeara_quadratic_1963},
	choose an adelic rotation $g=(g_{\mathfrak{p}})_{\mathfrak{p}\in\Omega_{F}}$ of $V$
	such that 
	\begin{align*}
		g_{\mathfrak{p}}L_{\mathfrak{p}}=L^{\prime}_{\mathfrak{p}}\quad\text{for}\quad\mathfrak{p}\not\in \infty_{F}\quad\text{and}\quad
		g_{\mathfrak{p}}=\id_{V_{\mathfrak{p}}}\quad\text{for}\quad \mathfrak{p}\notin \mathcal{T}.  
	\end{align*}

	For $\mathfrak{P}\mid\mathfrak{p}$, let $E_{\mathfrak{P}}$ be the completion of $E$
	at $\mathfrak{P}$ and write $\nr_{\mathfrak{P}\mid\mathfrak{p}}$ for the local norm
	\cite[\S15B]{omeara_quadratic_1963}. Let $\mathbb{I}_{E}$ be the idele group of $E$.
	Put $\sigma=g\perp\id_{\mathbb{H}}$. Then $\sigma M=M^{\prime}$.
	Identify $V$ with $(E,\gamma\nr_{E/F})$ for some $\gamma\in F^{\times}$.
	
	For each $\mathfrak{p}\in\mathcal{T}$, choose
	$\beta_{\mathfrak{P}}\in E_{\mathfrak{P}}^{\times}$ for $\mathfrak{P}\mid\mathfrak{p}$ as follows.
	If $\mathfrak{p}$ does not split in $E$, there is a unique $\mathfrak{P}\mid\mathfrak{p}$ and
	$V_{\mathfrak{p}}\cong(E_{\mathfrak{P}},\gamma\nr_{\mathfrak{P}\mid\mathfrak{p}})$.
	By \cite[55:2]{omeara_quadratic_1963}, there exist
	$x_{\mathfrak{p}},y_{\mathfrak{p}}\in E_{\mathfrak{P}}^{\times}$ such that
	\begin{align*}
		\theta(g_{\mathfrak{p}})
		&=Q(x_{\mathfrak{p}})Q(y_{\mathfrak{p}})F_{\mathfrak{p}}^{\times2}=\gamma^{2}\nr_{\mathfrak{P}\mid\mathfrak{p}}(x_{\mathfrak{p}})
		\nr_{\mathfrak{P}\mid\mathfrak{p}}(y_{\mathfrak{p}})F_{\mathfrak{p}}^{\times2}=\nr_{\mathfrak{P}\mid\mathfrak{p}}(\gamma x_{\mathfrak{p}}y_{\mathfrak{p}})
		F_{\mathfrak{p}}^{\times2}
	\end{align*}
	and put $\beta_{\mathfrak{P}}=\gamma x_{\mathfrak{p}}y_{\mathfrak{p}}$.
	If $\mathfrak{p}$ splits in $E$, let $\mathfrak{P},\mathfrak{P}^{\prime}$ be the two primes above it.
	Then $E_{\mathfrak{P}}=E_{\mathfrak{P}^{\prime}}=F_{\mathfrak{p}}$ and both local norm maps are the identity.
	Choose $u_{\mathfrak{p}}\in F_{\mathfrak{p}}^{\times}$ with
	$\theta(g_{\mathfrak{p}})=u_{\mathfrak{p}}F_{\mathfrak{p}}^{\times2}$, and put
	\begin{align*}
		\beta_{\mathfrak{P}}=u_{\mathfrak{p}}
		\quad\text{and}\quad
		\beta_{\mathfrak{P}^{\prime}}=1.
	\end{align*}
	For $\mathfrak{p}\notin\mathcal{T}$, put $\beta_{\mathfrak{P}}=1$ for every
	$\mathfrak{P}\mid\mathfrak{p}$. Since $\mathcal{T}$ is finite,
	$\beta=(\beta_{\mathfrak{P}})_{\mathfrak{P}\in\Omega_{E}}\in\mathbb{I}_{E}$.
	By \cite[55:4 and \S64]{omeara_quadratic_1963},
	\begin{align*}
		\theta(\sigma_{\mathfrak{p}})=\theta(g_{\mathfrak{p}})
		&=\left(\prod_{\mathfrak{P}\mid\mathfrak{p}}
		\nr_{\mathfrak{P}\mid\mathfrak{p}}(\beta_{\mathfrak{P}})\right)F_{\mathfrak{p}}^{\times2}=(\nr_{E/F}(\beta))_{\mathfrak{p}}F_{\mathfrak{p}}^{\times2}
	\end{align*}
	for $\mathfrak{p}\in\Omega_{F}$.
	Since $FM$ is isotropic, by \cite[55:2a]{omeara_quadratic_1963}, we have $
		\theta(O^{+}(FM))=F^{\times}/F^{\times2}$. 
	Since $E=H_{F}$, the extension $E/F$ is unramified at every finite prime
	and splits completely at every infinite prime. Hence, by
	\cite[63:16 and 65:2]{omeara_quadratic_1963}, we have
	\begin{align}\label{inclu}
		\mathbb{U}_{F}\mathbb{I}_{F}^{2}\subseteq\nr_{E/F}(\mathbb{I}_{E}).
	\end{align}
	Since $h(F)=2$, we have $\cl(F)^{2}=\{1\}$.
	By \eqref{clf2} and \cite[65:21]{omeara_quadratic_1963},
	we obtain 
	\begin{align}\label{index}
		[\mathbb{I}_{F}:F^{\times}\nr_{E/F}(\mathbb{I}_{E})]
		=2=h(F)=[\mathbb{I}_{F}:F^{\times}\mathbb{U}_{F}\mathbb{I}_{F}^{2}].
	\end{align}
	Combining \eqref{inclu}, \eqref{index} and Lemma \ref{lem:hyperbolic-spinor}(i), we see that
	\begin{align*}
		\nr_{E/F}(\beta)\in F^{\times}\nr_{E/F}(\mathbb{I}_{E})
		=F^{\times}\mathbb{U}_{F}\mathbb{I}_{F}^{2}\subseteq F^{\times}\mathbb{I}_{M}.
	\end{align*}
	Write $ \nr_{E/F}(\beta)=\alpha j$ with 
		$\alpha\in F^{\times}$ and
		$j=(j_{\mathfrak{p}})_{\mathfrak{p}\in\Omega_{F}}\in\mathbb{I}_{M}$. 
	Choose $\rho\in O^{+}(FM)$ and, for each 
	$\mathfrak{p}\in \Omega_{F}\setminus \infty_{F}$, $\tau_{\mathfrak{p}}\in O^{+}(M_{\mathfrak{p}})$ with
	\begin{align*}
		\theta(\rho)=\alpha F^{\times2}
		\quad\text{and}\quad
		\theta(\tau_{\mathfrak{p}})
		=j_{\mathfrak{p}}F_{\mathfrak{p}}^{\times2}.
	\end{align*}
	Then
	\begin{align*}
		\theta(\rho_{\mathfrak{p}}\tau_{\mathfrak{p}})
		&=\alpha j_{\mathfrak{p}}F_{\mathfrak{p}}^{\times2}
		=(\nr_{E/F}(\beta))_{\mathfrak{p}}F_{\mathfrak{p}}^{\times2}=\theta(\sigma_{\mathfrak{p}}).
	\end{align*}
	Put $\eta_{\mathfrak{p}}
	=\rho_{\mathfrak{p}}^{-1}\sigma_{\mathfrak{p}}\tau_{\mathfrak{p}}^{-1}$.
	Then
	\begin{align*}
		\theta(\eta_{\mathfrak{p}})=1
		\quad\text{and}\quad
		M^{\prime}_{\mathfrak{p}}
		=\rho_{\mathfrak{p}}\eta_{\mathfrak{p}}M_{\mathfrak{p}}.
	\end{align*}
	By the definition of the proper spinor genus and
	Lemma \ref{lem:hyperbolic-spinor}(ii),
	\begin{align*}
		M^{\prime}\in\spn^{+}(M)=\cls^{+}(M).
	\end{align*}
	Thus $L^{\prime}\perp\mathbf{H}\cong L\perp\mathbf{H}$.
\end{proof}

\begin{lem}\label{lem:residue-isotropy}
	Suppose that $L$ is integral and $n\ge1$.
	For each prime ideal $\mathfrak{p}$, let
	$\overline{Q}$ be the quadratic form induced by $Q$ on
	$L/\mathfrak{p}L$.
	There is a finite set $\mathcal{S}$ of prime ideals
	of $\mathcal{O}_{F}$ containing all dyadic primes such that
	$\overline{Q}$ is nonsingular over $\kappa_{\mathfrak{p}}$
	for every $\mathfrak{p}\notin\mathcal{S}$.
	Moreover, the following hold.
	\begin{enumerate}[label=\rm(\roman*), itemindent=0em, leftmargin=*]
		\item If $n\ge3$ and $\mathfrak{p}\notin\mathcal{S}$,
		then $\overline{Q}$ is isotropic.
		\item If $n=2$, $\mathfrak{p}\notin\mathcal{S}$, and
		$\mathfrak{p}$ splits in $E=F(\sqrt{-\det V})$,
		then $\overline{Q}$ is isotropic.
		\item If $n=1$ and $\mathfrak{p}$ divides
		$\mathfrak{n}(L)$, then $\overline{Q}=\overline{0}$
		and hence isotropic.
	\end{enumerate}
\end{lem}

\begin{proof}
	Since $L$ is integral and
	$2\mathfrak{s}(L)\subseteq\mathfrak{n}(L)$,
	\cite[82:10]{omeara_quadratic_1963} gives
	$2^{n+1}\mathfrak{v}(L)
	\subseteq2\mathfrak{n}(L)^{n}\subseteq2\mathcal{O}_{F}$.
	Let $\mathcal{S}$ be the finite set of prime ideals dividing
	$2^{n+1}\mathfrak{v}(L)$.  For $\mathfrak{p}\notin\mathcal{S}$,
	let $A_{\mathfrak{p}}$ be a Gram matrix of
	$L_{\mathfrak{p}}$. Since 
	$\mathfrak{p}\nmid2$,
	\begin{align*}
		A_{\mathfrak{p}}\in M_{n}(\mathcal{O}_{F_{\mathfrak{p}}})
		\quad\text{and}\quad
		(\det A_{\mathfrak{p}})\mathcal{O}_{F_{\mathfrak{p}}}
		=\mathfrak{v}(L)_{\mathfrak{p}}
		=\mathcal{O}_{F_{\mathfrak{p}}}.
	\end{align*}
	For every prime ideal $\mathfrak{p}$,
	$L_{\mathfrak{p}}$ is free of rank $n$ over
	$\mathcal{O}_{F_{\mathfrak{p}}}$. Hence
	\begin{align*}
		L/\mathfrak{p}L
		\cong L_{\mathfrak{p}}/\mathfrak{p}L_{\mathfrak{p}}
		\cong\kappa_{\mathfrak{p}}^{n}.
	\end{align*} 
	If $n=1$ and $\mathfrak{p}$ divides $\mathfrak{n}(L)$,
	then $Q(L)\subseteq\mathfrak{n}(L)\subseteq\mathfrak{p}$.
	Thus $\overline{Q}=0$ and every nonzero vector of
	$L/\mathfrak{p}L$ is isotropic. This proves (iii).
	
	For $\mathfrak{p}\notin\mathcal{S}$, the Gram matrix
	$\overline{A_{\mathfrak{p}}}$ is nonsingular, and $\kappa_{\mathfrak{p}}$
	has odd characteristic.
	If $n=2$ and $\mathfrak{p}$ splits in
	$E=F(\sqrt{-\det V})$, then
	\begin{align*}
		-\det A_{\mathfrak{p}}
		\in F_{\mathfrak{p}}^{\times2}
		\cap\mathcal{O}_{F_{\mathfrak{p}}}^{\times}
		=\mathcal{O}_{F_{\mathfrak{p}}}^{\times2}.
	\end{align*}
	Hence $
		-\det\overline{A_{\mathfrak{p}}}
		\in\kappa_{\mathfrak{p}}^{\times2}$.
	Assertions (i) and (ii) now follow from
	\cite[62:1b and 42:9]{omeara_quadratic_1963},
	respectively.
\end{proof}

\begin{lem}\label{lem:field-class-group}
	Suppose that $L$ is norm Euclidean of rank $n\ge1$,
	and put $E=F(\sqrt{-\det V})$.
	Then the following hold.
	\begin{enumerate}[label=\rm(\roman*), itemindent=0em, leftmargin=*]
		\item If $n\ge3$, then $H_{F}=F$.
		\item If $n=2$, then $H_{F}\subseteq E$.
		\item If $n=1$, write $L=\mathfrak{a}v$, where
		$v\in V\setminus\{0\}$ and $\mathfrak{a}$ is a nonzero
		fractional ideal of $F$. Then
		\begin{align*}
			\cl(F)=\langle[\mathfrak{a}]\rangle
			\quad\text{and}\quad
			[\mathfrak{a}]^{2}=1.
		\end{align*}
	\end{enumerate}
\end{lem}
\begin{proof}
	Let $\mathcal{S}$ be as in Lemma \ref{lem:residue-isotropy}.
	By Lemmas \ref{lem:residue-isotropy} and
	\ref{lem:isotropic-prime}, every prime ideal $\mathfrak{p}$ satisfying
	\begin{align*}
		\begin{cases}
			\mathfrak{p}\notin\mathcal{S},
			& \text{if $n\ge3$},\\
			\mathfrak{p}\notin\mathcal{S}
			\text{ and $\mathfrak{p}$ splits in $E$},
			& \text{if $n=2$},\\
			\mathfrak{p}\mid\mathfrak{n}(L),
			& \text{if $n=1$}
		\end{cases}
	\end{align*}
	is principal.
	
	To prove (iii), suppose that $n=1$ and write
	$L=\mathfrak{a}v$, where $v\in V\setminus\{0\}$ and
	$\mathfrak{a}$ is a nonzero fractional ideal of $F$.
	Put $\gamma=Q(v)\in F^{\times}$. Since $\mathfrak{a}$ is locally principal,
	\begin{align*}
		\mathfrak{b}:=\mathfrak{n}(L)=\gamma\mathfrak{a}^{2}
		\subseteq\mathcal{O}_{F}
		\quad\text{and}\quad
		\nr(\mathfrak{b})=|\nr_{F/\mathbb{Q}}(\gamma)|\nr(\mathfrak{a})^{2}\geq1.
	\end{align*}
	For every $x\in F$, by norm Euclideanity, there exists
	$y\in\mathfrak{a}$ such that
	\begin{align*}
		|\nr_{F/\mathbb{Q}}(Q((x-y)v))|=|\nr_{F/\mathbb{Q}}(\gamma(x-y)^{2})|<1.
	\end{align*}
	Hence
	\begin{align*}
		|\nr_{F/\mathbb{Q}}(x-y)|
		&<\dfrac{1}{\sqrt{|\nr_{F/\mathbb{Q}}(\gamma)|}}
		=\dfrac{\nr(\mathfrak{a})}{\sqrt{\nr(\mathfrak{b})}}
		\leq\nr(\mathfrak{a}).
	\end{align*}
	Thus $[\mathfrak{a}]$ is a norm Euclidean ideal class in the sense of
	\cite[(0.2)]{lenstra_euclidean_ideals_1979}. By \cite[(0.3) Theorem]{lenstra_euclidean_ideals_1979}, we have $
		\cl(F)=\langle[\mathfrak{a}]\rangle$.
	Since every prime ideal dividing
	$\mathfrak{b}=\mathfrak{n}(L)$ is principal, $\mathfrak{b}$ is principal and
	\begin{align*}
		[\mathfrak{a}]^{2}=[\gamma\mathfrak{a}^{2}]=[\mathfrak{b}]=1.
	\end{align*}
	This proves (iii).
	
	For (i) and (ii), assume that $n\ge2$ and put
	\begin{align*}
		K=
		\begin{cases}
			F,&\text{if $n\ge3$},\\
			E,&\text{if $n=2$}.
		\end{cases}
	\end{align*}
	Every prime ideal $\mathfrak{p}\notin\mathcal{S}$ that splits
	completely in $K$ is principal, and hence splits completely
	in $H_{F}$ by \cite[Corollary 5.25]{cox_primes_2013}.
	Since $K$ and $H_{F}$ are Galois over $F$,
	\cite[Theorem 8.19(i)]{cox_primes_2013} implies
	$H_{F}\subseteq K$. This proves (i) and (ii).
\end{proof}

\begin{proof}[Proof of Theorem \ref{thm:class-number}]
	To prove maximality, let $L^{\prime}$ be an integral
	$\mathcal{O}_{F}$-lattice on $V$ containing $L$.
	For any $x\in L^{\prime}$, there exists $y\in L$ such that $
		\left|\nr_{F/\mathbb{Q}}(Q(x-y))\right|<1$. Since $x-y\in L^{\prime}$ and $L^{\prime}$ is integral,
	$Q(x-y)\in\mathcal{O}_{F}$. Hence $
		\left|\nr_{F/\mathbb{Q}}(Q(x-y))\right|
		\in\mathbb{Z}_{\ge 0}$,
	so $Q(x-y)=0$. It follows from anisotropy that $x=y\in L$.
	Therefore, $L^{\prime}=L$ and thus $L$ is
	$\mathcal{O}_{F}$-maximal.

	By \cite[Theorem 5.23]{cox_primes_2013}, $\cl(F)\cong \gal(H_{F}/F)$ and thus
	$h(F)=[H_{F}:F]$.
	If $n\ge3$, by Lemma \ref{lem:field-class-group}(i),
	$H_{F}=F$, so $h(F)=1$.
	If $n=2$, since $V$ is anisotropic,
	$-\det V\notin F^{\times2}$, and hence $[E:F]=2$.
	By Lemma \ref{lem:field-class-group}(ii), $	F\subseteq H_{F}\subseteq E$. So
	\begin{align*}
	    [E:H_{F}]h(F)=[E:H_{F}][H_{F}:F]=[E:F]=2.
	\end{align*}
	Thus $h(F)\leq2$, and
	\begin{align*}
		h(F)=2
		\quad\Longleftrightarrow\quad [E:H_{F}]=1
		\quad\Longleftrightarrow\quad H_{F}=E.
	\end{align*}
	This proves (i) and (ii).
	If $n=1$, write $L=\mathfrak{a}v$ as in
	Lemma \ref{lem:field-class-group}(iii). Then
	$h(F)=|\langle[\mathfrak{a}]\rangle|\leq2$,
	with equality if and only if $\mathfrak{a}$ is not principal,
	or equivalently, if and only if $L$ is not free.
	This proves (iii). Since $h(F)\leq2$ in all cases,
	Lemma \ref{lem:stabilized-class-number} implies
	$h(F)=h(L\perp\mathbf{H})$.

	We now prove $h(L)=1$.
	If $n=1$, then $O(V_{\mathfrak{p}})=\{\pm\id_{V_{\mathfrak{p}}}\}$
	for every finite prime $\mathfrak{p}$, so every
	$L^{\prime}\in\gen(L)$ satisfies
	$L^{\prime}_{\mathfrak{p}}=L_{\mathfrak{p}}$ at every finite prime.
	Hence $L^{\prime}=L$ and $h(L)=1$.

	Assume that $n\ge2$. Let $L^{\prime}\in\gen(L)$. Put
	$M=L\perp\mathbf{H}$ and
	$M^{\prime}=L^{\prime}\perp\mathbf{H}\in\gen(M)$.
	If $h(F)=1$, then $h(M)=1$, so $M^{\prime}\cong M$.
	If $h(F)=2$, then by (i) and (ii), we have $n=2$ and $E=H_{F}$.
	Lemma \ref{lem:binary-stabilization} implies $M^{\prime}\cong M$.
	In either case, Theorem \ref{thm:cancellation} implies $L^{\prime}\cong L$.
	Therefore, $h(L)=1$.
\end{proof}
\begin{proof}[Proof of Corollary \ref{cor:stabilized-genus}]
	The proof of Lemma \ref{lem:stabilized-class-number} shows that
	$L\perp\mathbf{H}(\mathfrak{a})\in\gen(L\perp\mathbf{H})$
	and $L\perp\mathbf{H}(\mathfrak{a})\not\cong L\perp\mathbf{H}$. By the same lemma, $h(L\perp\mathbf{H})=h(F)=2$,
	so the assertion follows.
\end{proof}

\section{Finiteness results}\label{sec:finiteness}
Throughout this section, let $F$ be an algebraic number field
and $L$ an $\mathcal{O}_{F}$-lattice. For $\alpha\in F^{\times}$, write $L^{(\alpha)}$ for $L$
equipped with $\alpha Q$.
Since scaling preserves global and local isometries, we have
\begin{align*}
	h(L^{(\alpha)})=h(L);
\end{align*}
see \cite[\S82J]{omeara_quadratic_1963}.
\begin{lem}\label{lem:norm-bound}
	Let $F$ be an algebraic number field and $L$ be an integral norm Euclidean $\mathcal{O}_{F}$-lattice
	of rank $n\ge 1$.
	Then
	\begin{align*}
		\nr(\mathfrak{n}(L))<4^{d}.
	\end{align*} 
\end{lem}

\begin{proof}
	Since $L$ is a free $\mathbb{Z}$-module of rank $nd$, we have $L/2L\cong (\mathbb{Z}/2\mathbb{Z})^{nd}$, so
	$[L:2L]=2^{nd}>1$. Choose $x\in FL$ such that $2x\in L\setminus2L$. Since $L$ is norm Euclidean, there exists $y\in L$ such that $|\nr_{F/\mathbb{Q}}(Q(x-y))|<1$. Put
	$w=2x-2y \in L\setminus 2L$. Since $FL$ is anisotropic, $Q(w)\neq0$. Hence
	\begin{align*}
		0<|\nr_{F/\mathbb{Q}}(Q(w))|
		=4^{d}|\nr_{F/\mathbb{Q}}(Q(x-y))|
		<4^{d}.
	\end{align*}
	Since $Q(w)\mathcal{O}_{F}\subseteq\mathfrak{n}(L)$, we have
	\begin{align*}
		|\nr_{F/\mathbb{Q}}(Q(w))|
		&=[\mathcal{O}_{F}:Q(w)\mathcal{O}_{F}]\\
		&=[\mathcal{O}_{F}:\mathfrak{n}(L)]\,[\mathfrak{n}(L):Q(w)\mathcal{O}_{F}]=\nr(\mathfrak{n}(L))\,[\mathfrak{n}(L):Q(w)\mathcal{O}_{F}].
	\end{align*}
	Hence $\nr(\mathfrak{n}(L))\leq|\nr_{F/\mathbb{Q}}(Q(w))|<4^{d}$.
\end{proof}

\begin{re}\label{re:norm-bound-high-rank}
	In fact, if $L$ is an integral norm Euclidean $\mathcal{O}_{F}$-lattice
	of rank $n\ge3$, then
	\begin{align*}
		\mathfrak{n}(L)=\mathcal{O}_{F}
		\quad\text{and}\quad \nr(\mathfrak{n}(L))=1.
	\end{align*}
	Indeed, by Theorem \ref{thm:class-number}, we have $h(F)=1$, so
	$\mathcal{O}_{F}$ is a PID and $L$ is free. Expressed in an
	$\mathcal{O}_{F}$-basis of $L$, the quadratic form $Q$ is ADC
	by \cite[Theorem 8]{clark_ADC-I-2012}, that is,
	$Q(V)\cap\mathcal{O}_{F}=Q(L)$.
	Since $\mathcal{O}_{F}$ is a Hasse domain and $n\ge3$, this form is
	primitive by \cite[Theorem 3.5]{clark_ADC-II-2014}, which means
	$\mathfrak{n}(L)=\mathcal{O}_{F}$.
\end{re}

We reformulate Pfeuffer's finiteness theorems in terms of the norm ideal.
\begin{prop}\label{prop:pfeuffer-norm}
	Let $h_{0}$ be a positive integer.
	\begin{enumerate}[label=\rm(\roman*), itemindent=0em, leftmargin=*]
		\item There are only finitely many isomorphism classes of totally
		real fields $F$ admitting a positive definite integral
		$\mathcal{O}_{F}$-lattice $L$ of rank $n\ge3$ with $h(L)\leq h_{0}$.
		\item For every positive integer $d_{0}$, there are only
		finitely many isomorphism classes of totally real fields $F$
		of degree at most $d_{0}$ admitting a positive definite
		integral $\mathcal{O}_{F}$-lattice $L$ of rank $n=2$ with
		$h(L)\leq h_{0}$.
		\item For a fixed totally real field $F$ and a fixed $C>0$,
		there are only finitely many isometry classes of positive
		definite integral $\mathcal{O}_{F}$-lattices $L$ of arbitrary
		rank $n\ge2$ satisfying
		\begin{align*}
			h(L)\leq h_{0}
			\quad\text{and}\quad
			\nr(\mathfrak{n}(L))\leq C.
		\end{align*}
	\end{enumerate}
\end{prop}
\begin{proof}
	Recall from \cite[\S82E]{omeara_quadratic_1963} that $
		2\mathfrak{s}(L)\subseteq\mathfrak{n}(L)
		\subseteq\mathfrak{s}(L)$. Since $L$ is integral,
	\begin{align*}
		\mathfrak{s}(L^{(2)})=2\mathfrak{s}(L)
		\subseteq\mathfrak{n}(L)\subseteq\mathcal{O}_{F}.
	\end{align*}
	 
	 Thus $L^{(2)}$ is a positive definite lattice with integral scale. Also, $h(L^{(2)})=h(L)\leq h_{0}$.

	(i) Applying \cite[Theorem 1]{pfeuffer_class-numbers_1979}
	to $L^{(2)}$, we obtain the assertion.
	
	(ii) Applying \cite[Theorem 2]{pfeuffer_class-numbers_1979}
	to $L^{(2)}$ in each degree at most $d_{0}$, we obtain the assertion.
	
	(iii) Put $d=[F:\mathbb{Q}]$. Then
	\begin{align*}
		\nr(\mathfrak{s}(L^{(2)}))
		&=\nr(2\mathfrak{s}(L))
		\leq\nr(2\mathfrak{n}(L))=2^{d}\nr(\mathfrak{n}(L))
		\leq2^{d}C.
	\end{align*}
	Hence $\mathfrak{s}(L^{(2)})$ ranges over finitely many
	integral ideals
	\cite[Theorem 12.5.3]{alaca_introductory_2004}.
	For each such scale, \cite[Theorem 1]{pfeuffer_class-numbers_1979}
	in ranks $n\ge3$ and \cite[Theorem 3]{pfeuffer_class-numbers_1979}
	in rank $2$ give only finitely many isometry classes of $L^{(2)}$.
	Since $L^{(2)}\cong L^{\prime (2)}$ if and only if $L\cong L^{\prime}$,
	the same finiteness holds for $L$.
\end{proof}
 
The following lemma is the fixed-field case of
\cite[Lemma 5.1]{he_trace_2026}. We include a short proof for completeness.
\begin{lem}\label{lem:bounded-volume-finiteness}
	Fix a totally real field $F$, a positive integer $n$, and a real
	number $C>0$. There are only finitely many isometry classes of
	positive definite integral $\mathcal{O}_{F}$-lattices $L$ of rank $n$
	satisfying
	\begin{align*}
		\nr(\mathfrak{v}(L))\leq C.
	\end{align*}
\end{lem}
\begin{proof}
	Put $d=[F:\mathbb{Q}]$. By \cite[\S 82E and \S 82J]{omeara_quadratic_1963},
	\begin{align*}
		\mathfrak{s}(L^{(2)})
		&=2\mathfrak{s}(L)\subseteq\mathfrak{n}(L)
		\subseteq\mathcal{O}_{F}\quad\text{and}\quad 
		\mathfrak{v}(L^{(2)})
		=2^{n}\mathfrak{v}(L)
		\subseteq\mathfrak{s}(L^{(2)})^{n}\subseteq\mathcal{O}_{F}.
	\end{align*}
	Hence $L^{(2)}$ has integral scale and integral volume, with
	\begin{align*}
		\nr(\mathfrak{v}(L^{(2)}))
		=2^{nd}\nr(\mathfrak{v}(L))\leq2^{nd}C.
	\end{align*}
	Thus $\mathfrak{v}(L^{(2)})$ ranges over finitely many integral
	ideals \cite[Theorem 12.5.3]{alaca_introductory_2004}.
	Applying \cite[Remark 103:5 and Theorem 103:4]{omeara_quadratic_1963}
	successively, we obtain only finitely many isometry classes of $L^{(2)}$.
	Since $L^{(2)}\cong L^{\prime (2)}$ if and only if $L\cong L^{\prime}$,
	the same holds for $L$.
\end{proof}

\begin{proof}[Proof of Theorem \ref{thm:finiteness}]
	For (i) and (ii), let $L$ be a lattice under consideration. By
	Theorem \ref{thm:class-number}, we have $h(L)=1$. 
	Choose $v\in L\setminus\{0\}$ and put $\alpha=Q(v)$.
	Since $L$ is totally definite and integral,
	$\alpha\in\mathcal{O}_{F}\setminus\{0\}$ has the same sign
	as $Q$ at every real embedding.
	Thus $L^{(\alpha)}$ is positive definite and integral,
	with $h(L^{(\alpha)})=1$.
	Applying Proposition \ref{prop:pfeuffer-norm}(i) with $h_{0}=1$ and Proposition \ref{prop:pfeuffer-norm}(ii) with $(h_{0},d_{0})=(1,d)$, we see that only finitely many fields $F$ occur,
	up to isomorphism.
	
	It remains to prove finiteness over a fixed totally real field $F$.
	Put $d=[F:\mathbb{Q}]$.
	For each sign pattern $\varepsilon$ that occurs, fix
	$\alpha_{\varepsilon}\in\mathcal{O}_{F}\setminus\{0\}$
	having that sign pattern, as above.
	For every lattice $L$ under consideration with sign pattern
	$\varepsilon$, put $J=L^{(\alpha_{\varepsilon})}$.
	Then $J$ is positive definite and integral, with
	$h(J)=h(L)=1$.
	By Lemma \ref{lem:norm-bound},
	\begin{align*}
		\nr(\mathfrak{n}(J))
		&=|\nr_{F/\mathbb{Q}}(\alpha_{\varepsilon})|\,
		\nr(\mathfrak{n}(L)) <4^{d}|\nr_{F/\mathbb{Q}}(\alpha_{\varepsilon})|
		=:C_{\varepsilon}.
	\end{align*}
	If $\rank\,L\ge2$, by Proposition \ref{prop:pfeuffer-norm}(iii)
	with $(h_{0},C)=(1,C_{\varepsilon})$, there are only finitely many
	isometry classes of $J$.
	If $\rank\,L=1$, then $\mathfrak{v}(J)=\mathfrak{n}(J)$,
	so by Lemma \ref{lem:bounded-volume-finiteness} with
	$(n,C)=(1,C_{\varepsilon})$, we obtain the same conclusion.
	Since $\alpha_{\varepsilon}$ is fixed, the same holds for $L$.
	There are at most $2^{d}$ sign patterns, so there are only finitely many
	isometry classes over each fixed $F$.
	This proves (iii) and, by the preceding discussion, (i) and (ii).
\end{proof}

\begin{proof}[Proof of Theorem \ref{thm:rank-degree-bounds}]
	By Theorem \ref{thm:class-number}, $L$ is a positive definite
	$\mathcal{O}_{F}$-maximal lattice of rank $n\ge3$ and class number one.
	Moreover, $L$ is primitive by Remark \ref{re:norm-bound-high-rank}.
	If $F=\mathbb{Q}$, then $[F:\mathbb{Q}]=1$ and $n\leq10$
	by Hanke's classification \cite[Theorem 7.5 and \S8]{hanke_maximal_2011}.
	If $F\neq\mathbb{Q}$, then $n\leq6$ and $[F:\mathbb{Q}]\leq5$
	by Kirschmer's classification \cite[Theorem 7.1 and \S7]{kirschmer_one-class_2014}.
\end{proof} 
\section*{Acknowledgments}
This work was supported by the National Natural Science Foundation of China (Project No. 12301013 and No. 12671016).

\section*{AI Declaration}
ChatGPT (GPT-5.6 Sol) was used to assist with exploring the proof strategy, using an earlier version of the manuscript as a starting point. It was also used to proofread the manuscript, suggest improvements to the exposition, and develop the Lean 4 formalization of selected main results \cite{he_norm-euclidean-verification_2026}. The author identified the mathematical connections and organized and wrote the manuscript. All mathematical statements and proofs were reviewed and verified by the author, who takes full responsibility for the content of this paper.

\begin{sortedbibliography}{99}
\bibitem[Alaca Williams 2004]{alaca_introductory_2004}
S.~Alaca and K.~S. Williams, \textit{Introductory Algebraic Number Theory},
Cambridge University Press, Cambridge, 2004. 

\bibitem[Borcherds 1987]{borcherds_lorentzian_1987}
R.~E. Borcherds, \textit{Automorphism groups of Lorentzian lattices},
J. Algebra \textbf{111} (1987), 133--153.

\bibitem[Clark 2012]{clark_ADC-I-2012}
P.~L. Clark, \textit{Euclidean quadratic forms and ADC forms: I},
Acta Arith. \textbf{154} (2012), 137--159. 

\bibitem[Clark Jagy 2014]{clark_ADC-II-2014}
P.~L. Clark and W.~C. Jagy,
\textit{Euclidean quadratic forms and ADC-forms II: integral forms},
Acta Arith. \textbf{164} (2014), 265--308. 

\bibitem[Cox 2013]{cox_primes_2013}
D.~A. Cox, \textit{Primes of the Form $x^{2}+ny^{2}$: Fermat, Class Field Theory,
and Complex Multiplication}, 2nd ed., John Wiley \& Sons, Hoboken, NJ, 2013. 

\bibitem[Eichler 1974]{eichler_quadratische_1974}
M. Eichler, \textit{Quadratische Formen und orthogonale Gruppen},
2nd ed., Grundlehren Math. Wiss. 63, Springer, Berlin, 1974.

\bibitem[Hanke 2011]{hanke_maximal_2011}
J. Hanke, \textit{Enumerating maximal definite quadratic forms of
bounded class number over $\mathbb{Z}$ in $n\ge3$ variables},
(2011), arXiv:1110.1876.

\bibitem[He 2026]{he_trace_2026}
Z. He,
\textit{Finiteness theorems for trace Euclidean quadratic lattices over totally real number fields}, (2026), arXiv:2610.34604.

\bibitem[He 2026 Verification]{he_norm-euclidean-verification_2026}
Z. He,
\textit{Norm-Euclidean verification package},
GitHub repository, 2026,
\url{https://github.com/hzldew-git/PaperVerification/tree/main/Norm-Euclidean}.

\bibitem[He Hu Xu 2023]{he-hu-xu_indefinite_2023}
Z. He, Y. Hu and F. Xu,
\textit{On indefinite $k$-universal integral quadratic forms over
number fields}, Math. Z. \textbf{304} (2023), Paper No. 20, 26 pp. 

\bibitem[Kirschmer 2014]{kirschmer_one-class_2014}
M. Kirschmer, \textit{One-class genera of maximal integral quadratic
forms}, J. Number Theory \textbf{136} (2014), 375--393.

\bibitem[Lenstra 1979]{lenstra_euclidean_ideals_1979}
H.~W. Lenstra, Jr., \textit{Euclidean ideal classes},
Ast\'erisque \textbf{61} (1979), 121--131.

\bibitem[OMeara 2000]{omeara_quadratic_1963}
O.~T. O'Meara, \textit{Introduction to Quadratic Forms},
Classics in Mathematics, Springer, Berlin, 2000.

\bibitem[OMeara Pollak 1965]{omeara_generation_1965}
O.~T. O'Meara and B. Pollak,
\textit{Generation of local integral orthogonal groups},
Math. Z. \textbf{87} (1965), 385--400. 

\bibitem[Pfeuffer 1979]{pfeuffer_class-numbers_1979}
H. Pfeuffer, \textit{On a conjecture about class numbers of totally
positive quadratic forms in totally real algebraic number fields},
J. Number Theory \textbf{11} (1979), 188--196.

\bibitem[Schulze-Pillot 1982]{schulze-pillot_remark_1982}
R. Schulze-Pillot, \textit{Remark on a note of Hsia about classification
of unimodular lattices}, J. Number Theory \textbf{14} (1982), 83--85.

\bibitem[Xu Zhang 2022]{xu-zhang_indefinite_2022}
F. Xu and Y. Zhang,
\textit{On indefinite and potentially universal quadratic forms over
number fields}, Trans. Amer. Math. Soc. \textbf{375} (2022), 2459--2480. 
\end{sortedbibliography}
\end{document}